\documentclass[12pt]{article}

\usepackage{amsmath}
\usepackage{amssymb}
\usepackage{amscd}

\newcommand{\bA}{{\bf A}}
\newcommand{\bB}{{\bf B}}
\newcommand{\bC}{{\bf C}}

\newcommand{\bP}{{\bf P}}
\newcommand{\bQ}{{\bf Q}}
\newcommand{\bR}{{\bf R}}
\newcommand{\bX}{{\bf X}}
\newcommand{\bZ}{{\bf Z}}

\newcommand{\CR}{\nonumber \\}
\def\br#1#2{\bX_{\bf [#1,#2]}}
\def\u1#1{{\langle #1 \rangle}}
\newtheorem{Cor}{Corollary}[section]
\newtheorem{Exa}{Example}[section]
\newtheorem{Lem}{Lemma}[section]
\newtheorem{Pro}{Proposition}[section]
\newtheorem{Rem}{Remark}[section]
\newtheorem{The}{Theorem}[section]
\newcommand{\ake}{\vskip 3mm \noindent}
\newcommand{\aki}{\vskip 5mm \noindent}

\begin{document}

\begin{titlepage}
\begin{flushright}

\end{flushright}
\vspace{0.5cm}
\begin{center}
{\Large \bf 
 Monodromies of rational elliptic surfaces and extremal elliptic $K3$ surfaces 
\par}
\lineskip .75em
\vskip2.5cm
\normalsize
{\large Mitsuaki Fukae}
\vskip 1.5em
{\large\it Department of Mathematics, Kobe University \\
Rokko, Kobe 657-8501, Japan}
\end{center}
\vskip3cm
\begin{abstract}
We give a systematic method to calculate
some homological data from the global monodromy of
a topological elliptic surface.
We apply this method to the cases 1) the transcendental 
lattice of an extremal elliptic $K3$ surface, 2) the torsion part of 
Mordell-Weil group of a general elliptic surface,
and 3) the Mordell-Weil lattice of a rational elliptic surface.

We present possible list of global monodromies of 
rational elliptic surface and extremal elliptic $K3$ surface
which reproduce to the classification obtained by
Oguiso-Shioda and Shimada-Zhang respectively.
\end{abstract}
\end{titlepage}

\tableofcontents

\baselineskip=0.7cm
\section{Introduction}
In our previous paper \cite{FYY} we analyze the trivial 
lattices (which is essentially the combinations of the singular fibres), 
the Mordell-Weil lattices and the torsion parts of Mordell-Weil group 
of rational elliptic surfaces using the method of string junctions in 
string theory, and reproduce the classification of the Mordell-Weil lattices 
due to Oguiso-Shioda by constructing the corresponding global monodromy. 
In this paper we apply the idea of \cite{FYY} to elliptic $K3$ surfaces. \\
\indent Let $f:X \rightarrow \bP^1$ be a general topological elliptic surface 
which has a smooth cross section and no multiple fibre. 
Let $T$ be the sublattice of $H_2(X,\bZ)$ which is generated by 
irreducible components of each singular fibre, a smooth cross section and 
a generic fibre. 
We denote the orthogonal complement of $T$ in $H_2(X,\bZ)$ and 
the primitive closure of $T$ in $H_2(X,\bZ)$ by $T^{\perp}$, $\hat{T}$ 
respectively. \\
The results of this paper is as follows : 
\begin{enumerate}
\item We give a systematic method to calculate the homological data $T^{\perp}$, 
$\hat{T} / T$ from the global monodromies. 
\item We reproduce the classification of rational elliptic surfaces and 
      extremal elliptic $K3$ surfaces due to Oguiso-Shioda and Shimada-Zhang 
      by constructing the corresponding global monodromy. 
\end{enumerate}
\underline{About $1.$}\, 
We give these algorithms in section $\ref{algo}$. 
These algorithms are essentially reduced to the calculation of the 
elementary divisor of integral matrix. 
\ake
\underline{About $2.$}\, 
Let $f:X \rightarrow \bP^1$ be a complex elliptic $K3$ surface with a 
section $O$ which satisfies the following conditions :\\
\,\,\,\,\,\,(i)\, The Picard number is $20$. \\
\,\,\,\,\,(ii)\, The Mordell-Weil group is finite. \\
Then $f:X \rightarrow \bP^1$ is said to be an \textit{extremal} elliptic 
$K3$ surface. 

The classification of the pairs (trivial lattice, Mordell-Weil group) 
of semi-stable extremal elliptic $K3$ surface has been done by Miranda 
and Persson \cite{MP} and complemented by Artal-Bartolo, Tokunaga and 
Zhang \cite{ABTZ}. 
The classification of the pairs (trivial lattice, Mordell-Weil group) 
of unsemi-stable extremal elliptic $K3$ surface has been done by Ye \cite{Y}. 
The complete classification of the triplets (trivial lattice, 
transcendental lattice, Mordell-Weil group) of extremal elliptic $K3$ 
surfaces has been done by Shimada and Zhang \cite{SZ}. 

One of the purpose of this paper is to give the list of global monodromies 
which corresponds to the list given by Shimada and Zhang (called \textit
{Shimada-Zhang's table}). 
The main theorem of this paper is as follows : \\
\textbf{Theorem.}\,
If ther exist an extremal elliptic $K3$ surface with the global monodromy 
given in Table $\ref{config2}$, then the corresponding triplet 
(trivial lattice, transcendental lattice, Mordell-Weil group) of the 
surface is given by the data of the same Table $\ref{config2}$

\section{Preliminaries}

\subsection{Monodromy of elliptic surfaces}
In this section we will recall some basic facts on the monodromy of 
elliptic surfaces following the reference \cite{M}. 

Let $f:X \to \bP^{1}$ be an elliptic surface that has $n$ singular fibres 
$F_{u_1},F_{u_2},\dots,F_{u_n}$ on $S=\left\{u_1,u_2,\dots,u_n\right\}$.
We take an base point $u_0 \in {\bP}^{1} \setminus S$, and 
let $\rho_{1},\rho_{2},\dots,\rho_{n}$ be disjoint smooth paths which 
connect $u_0$ with $u_1,u_2,\dots,u_n$ respectively, 
$\rho'_{1},\rho'_{2},\dots,\rho'_{n}$ be generators 
of $\pi_{1}\left({\bP}^{1} \setminus S,u_0\right)$ which correspond to 
$\rho_{1},\rho_{2},\dots,\rho_{n}$,
$\rho_{1}^{*},\rho_{2}^{*},\dots,\rho_{n}^{*}$ be automorphisms of 
${H}_{1}\left(f^{-1}(u_0),\bZ \right)$ corresponding to 
$\rho'_{1},\rho'_{2},\dots,\rho'_{n}$ and 
$\langle \alpha,\beta \rangle$ be the basis of 
${H}_{1}\left(f^{-1}(u_0),\bZ\right)$.
We put $K_i\left(1\le i\le n\right)$ the circuit matrices along the 
path ${\rho'}_{i}$, and call $\left(K_1,K_2,\dots,K_n\right)$ the 
\textit{global monodromy} on $\left(\rho'_{1},\rho'_{2},\dots,\rho'_{n}\right)$ 
of $X$. 
Namely if ${\rho}_{i}^{*}\left(\alpha\right) = a_i\alpha+c_i\beta$ and 
${\rho}_{i}^{*}\left(\beta\right) = b_i\alpha+d_i\beta$ then 
$K_i = \begin{pmatrix} a_i & b_i \\ c_i & d_i \end{pmatrix}$.\footnote{We remark 
that this definition of the circuit matrix is inverse of the ordinary one.} 
It is clear also that 
$K_n K_{n-1}\dotsb K_1 = I.$
For $i=1,2,\dots,n-1$ let $u_i:\pi_{1}({\bP}^{1} \setminus S,u_0) \to 
\pi_{1}\left({\bP}^{1} \setminus S,u_0\right)$ be an automorphism of 
$\pi_{1}\left({\bP}^{1} \setminus S,u_0\right)$ defined by 
$u_i\left({\rho'}_{1}\right)={\rho'}_{1},\dots,
u_i\left({\rho'}_{i-1}\right)={\rho'}_{i-1},
u_i\left({\rho'}_{i}\right)={\rho'}_{i+1},
u_i\left({\rho'}_{i+1}\right)=
{\rho'}_{i+1}^{-1}\cdot{\rho'}_{i}\cdot{\rho'}_{i+1},
u_i\left({\rho'}_{i+2}\right)={\rho'}_{i+2},\dots,
u_i\left({\rho'}_{n}\right)={\rho'}_{n}$ 
and $R_i\left(R_i^{-1}\right)$ be an transformation of global monodromy 
which is induced by $u_i\left(u_i^{-1}\right)$. Then we have 
\begin{equation}
\begin{split}
& R_i\left(\left(K_1,\dots,K_n\right)\right)
 =\left(K_1,\dots,K_{i-1},K_{i+1},K_{i+1}K_{i}K_{i+1}^{-1},K_{i+2},
 \dots,K_n\right) 
 \\
& R_i^{-1}\left(\left(K_1,\dots,K_n\right)\right)
 =\left(K_1,\dots,K_{i-1},K_{i}^{-1}K_{i+1}K_{i},K_{i},K_{i+2},
 \dots,K_n\right). \label{elem}
\end{split}
\end{equation}
We call the transformations $R_i,R_i^{-1}\left(1\le i\le n-1\right)$ the
\textit{elementary transformations} of global monodromy. 
If there exists a finite sequence of elementary transformations 
and $SL(2,\bZ)$ transformations corresponding to the base change of 
$H_1(f^{-1}(u_0),\bZ)$ which 
transfer $\left(K_1,K_2,\dots,K_n\right)$ to 
$(K'_1,K'_2,\dots,K'_n)$, then we identify 
$\left(K_1,K_2,\dots,K_n\right)$ with 
$(K'_1,K'_2,\dots,K'_n)$
and denote 
$\left(K_1,K_2,\dots,K_n\right)\sim (K'_1,K'_2,
\dots,K'_n)$. 

\subsection{Oguiso-Shioda's table and Shimada-Zhang's table}
In this section we recall the results of the references \cite{OS} and 
\cite{SZ}.\\
Let $f:X_k \to {\bP}^{1}$ be a topological elliptic surface that has 
$n$ singular fibres 
$F_{u_1},F_{u_2},\dots,F_{u_n}$ on $S=\left\{u_1,u_2,\dots,u_n\right\}$ 
and the Euler number is $12k$. 
Through out of this paper we assume $f$ has a smooth cross section $O$. 
Let $O$, $F$ be a section and a generic fibre. 
We put $U' = \langle O,F \rangle$. 
Then it is known that 
\begin{equation*}
 H_2\left(X_k,\bZ\right) \cong \left(-E_8\right)^{\oplus k} \oplus 
U^{\oplus 2k-2}\oplus U',
\end{equation*}
where $U = \begin{pmatrix} 0 & 1 \\ 1 & 0\end{pmatrix}$ and  
$U' = \begin{pmatrix} -k & 1 \\ 1 & 0\end{pmatrix}.$
We set the sublattice $L_k$ of $H_2\left(X_k,\bZ\right)$;
\begin{equation*}
 L_k := \left(-E_8\right)^{\oplus k} \oplus U^{\oplus 2k-2}. 
\end{equation*}
For each $v \in S=\left\{u_1,u_2,\dots,u_n\right\}$, let 
\begin{equation*}
 F_v = f^{-1}\left(v\right) = \Theta_{v,0} + 
       \sum_{i=1}^{m_v-1}\mu_{v,i}\Theta_{v,i},
\end{equation*}
where $\Theta_{v,i}\left(0\le i\le m_v-1\right)$ are the irreducible 
components of $F_v$, $m_v$ being their multiplicity. 
$\Theta_{v,0}$ is the unique component of $F_v$ meeting the zero section. 
And we set the sublattice of $L_k$
\begin{equation*}
 V := \bigoplus_{v \in S} T_v \>,
\end{equation*}
where $T_v = \langle \Theta_{v,i} \mid 1\le i \le m_v - 1\rangle.$
Then it is known that the opposite lattice $-T_{v}$ is a root lattice 
of rank $m_v - 1$ which is determined by the type of the singular fibre 
$F_v$ \cite{Kod}. 
We can summarize the correspondence of the type of $F_v$ and $T_v$ in 
Table $\ref{table1}$. 
\begin{table}
\begin{center}
\begin{tabular}{|c|c|} \hline
\multicolumn{1}{|c|}{the type of $F_v$ } & $T_v$ \\ \hline
$I_n\left(1\le n\right)$ & $-A_{n-1}\left(A_0 = 0\right)$ \\ \hline
\setcounter{part}{2} $\thepart$ & $0$ \\ \hline
\setcounter{part}{3} $\thepart$ & $-A_1$ \\ \hline
\setcounter{part}{4} $\thepart$ & $-A_2$ \\ \hline
$I_{n}^{*}\left(0\le n\right)$ & $-D_{n+4}$ \\ \hline
\setcounter{part}{2} ${\thepart}^{*}$ & $-E_8$ \\ \hline
\setcounter{part}{3} ${\thepart}^{*}$ & $-E_7$ \\ \hline
\setcounter{part}{4} ${\thepart}^{*}$ & $-E_6$ \\ \hline
\end{tabular} 
\end{center}
\caption{the correspondence of the type of $F_v$ and $T_v$}
\label{table1}
\end{table}
\ake
\ake
\ake
Consider the lattice embedding 
\begin{equation}
i_V : V \hookrightarrow L_k \>, \label{emb1}
\end{equation}
and let $V^\perp$ be the orthogonal complement of $V$ in $L_k$ and $\hat{V}$ 
the primitive closure of $V$ in $L_k$.
Now we describe the geometric meaning of 
$V^{\perp}$, and $\hat{V}/V$. 
First, with respect to the geometric meaning of $\hat{V}/V$ the 
following is known. 
\ake
\textbf{Theorem 1.}(\cite{Shio})
The torsion part of Mordell-Weil group of $X_k$ is isomorphic to $\hat{V}/V$.
\ake
\ake
\ake
\ake
For $k=1$, $X_k$ is a rational elliptic surface, and the following is 
known. 
\ake
\textbf{Theorem 2.}(\cite{OS})
The Mordell-Weil lattice of a rational elliptic surface is 
isomorphic to the opposite lattice of ${\left(V^{\perp}\right)}^{*}$. 
\ake
One of the main results of the reference \cite{OS} is the classification 
of all the possible triplets 
\[
 \left(V,\>\textrm{Mordell-Weil lattice},\>\textrm{torsion part of 
 Mordell-Weil group}\right) 
\]
of rational elliptic surfaces. 
(We call this table \textit{Oguiso-Shioda's table}.)
\ake
\ake
\ake
\ake
For $k=2$, $X_k$ is an elliptic $K3$ surface. 
In particular an elliptic $K3$ surface is said to be \textit{extremal} 
if the Picard number is $20$ and Mordell-Weil group is finite. 
For an extremal elliptic $K3$ surface it is known that  
\ake
\textbf{Theorem 3.}(\cite{SZ})
The transcendental lattice of an extremal elliptic $K3$ surface is 
isomorphic to $V^{\perp}$. 
And the Mordell-Weil group is isomorphic to $\hat{V}/V$.
\ake
One of the main results of the reference \cite{SZ} is to classify all 
the possible triplets 
\[
\left(V,\>\textrm{transcendental lattice},\>\textrm{Mordell-Weil group}\right) 
\]
of extremal elliptic $K3$ surfaces. 
(We call this table \textit{Shimada-Chang's table}.) 

\section{Junction}

\subsection{Tadpole junction lattice and rational tadpole junction lattice}
In this section we define tadpole junction lattice and rational tadpole 
junction lattice by 
using the notations of the previous section and the descriptions of the 
reference \cite{MSY}.

Let $D_i\left(1\le i\le n\right)$ be small closed disjoint $2$-disks in 
$\bP^{1}$ with the center $u_i\left(1\le i\le n\right)$. 
We can assume that the sets 
$\rho_{i} \cap D_j\left(1\le i\le n,1\le j\le n,i\neq j\right)$ are 
empty and that each of the sets 
$\rho_{i} \cap \partial D_i\left(1\le i\le n\right)$ has only one point. 
Let $u'_i=\rho_{i} \cap \partial D_i$ and $\bar{\rho}_{i}$ be the 
part of $\rho_{i}$ from $u_0$ to $u'_i$.
We can assume that 
$\rho'_{i} = \bar{\rho}_{i} \cdot \partial D_i \cdot 
\bar{\rho}_{i}^{-1}$ 
($\partial D_i$ is oriented as the boundary of $D_i$).

Now we construct $2$-cycles on $X$. For a curve $\rho$ starting 
{}from $u_0$ and an $1$-cycle $\gamma$ on $f^{-1}(u_0)$, by the locally 
flatness of the fibering $f:X \to \bP^{1}$, there is a continuous 
family of $1$-cycle $\gamma_u$ of $f^{-1}\left(u\right)\left(u \in \rho \right)$ 
such that $\gamma_{u_0} = \gamma$; the union $\bigcup_{u \in \rho}\gamma_u$, 
which is a real 2-dimensional surface, is denoted by 
$\rho\left(\gamma\right)$. It is clear that 
\begin{equation}
 \rho\mu\left(\gamma\right) = \rho\left(\gamma\right) 
 + \mu\left(\rho^{*}\gamma\right)\quad
\left(\textrm{as the element of}\>H_2(X,\bR)\right) \label{rho}
\end{equation}

\begin{figure}[htbp]
  \begin{center}
\unitlength 0.1in
\begin{picture}(51.82,27.36)(1.60,-31.46)
%
\special{pn 8}%
\special{ar 1022 2018 480 480  0.0000000 6.2831853}%
%
\special{pn 8}%
\special{ar 2942 2002 480 480  0.0000000 6.2831853}%
%
\special{pn 8}%
\special{ar 4862 2002 480 480  0.0000000 6.2831853}%
%
\special{pn 8}%
\special{pa 2942 410}%
\special{pa 2942 1522}%
\special{fp}%
\special{sh 1}%
\special{pa 2942 1522}%
\special{pa 2962 1455}%
\special{pa 2942 1469}%
\special{pa 2922 1455}%
\special{pa 2942 1522}%
\special{fp}%
\special{pa 2934 418}%
\special{pa 4862 1530}%
\special{fp}%
\special{sh 1}%
\special{pa 4862 1530}%
\special{pa 4814 1479}%
\special{pa 4816 1503}%
\special{pa 4794 1514}%
\special{pa 4862 1530}%
\special{fp}%
\special{pa 2942 418}%
\special{pa 1022 1530}%
\special{fp}%
\special{sh 1}%
\special{pa 1022 1530}%
\special{pa 1090 1514}%
\special{pa 1068 1503}%
\special{pa 1070 1479}%
\special{pa 1022 1530}%
\special{fp}%
%
\special{pn 8}%
\special{pa 1014 2002}%
\special{pa 1014 3130}%
\special{dt 0.045}%
\special{pa 1014 3130}%
\special{pa 1014 3129}%
\special{dt 0.045}%
\special{pa 2934 2010}%
\special{pa 2934 3146}%
\special{dt 0.045}%
\special{pa 2934 3146}%
\special{pa 2934 3145}%
\special{dt 0.045}%
\special{pa 4854 2002}%
\special{pa 4854 3122}%
\special{dt 0.045}%
\special{pa 4854 3122}%
\special{pa 4854 3121}%
\special{dt 0.045}%
%
\special{pn 8}%
\special{pa 1582 2002}%
\special{pa 1694 2002}%
\special{fp}%
%
\special{pn 8}%
\special{pa 1766 2002}%
\special{pa 1878 2002}%
\special{fp}%
%
\special{pn 8}%
\special{pa 1966 2002}%
\special{pa 2078 2002}%
\special{fp}%
%
\special{pn 8}%
\special{pa 2182 2002}%
\special{pa 2294 2002}%
\special{fp}%
%
\special{pn 8}%
\special{pa 2350 2002}%
\special{pa 2462 2002}%
\special{fp}%
%
\special{pn 8}%
\special{pa 3478 2010}%
\special{pa 3590 2010}%
\special{fp}%
%
\special{pn 8}%
\special{pa 3702 2010}%
\special{pa 3814 2010}%
\special{fp}%
%
\special{pn 8}%
\special{pa 3950 2010}%
\special{pa 4062 2010}%
\special{fp}%
%
\special{pn 8}%
\special{pa 4166 2010}%
\special{pa 4278 2010}%
\special{fp}%
\put(21.0000,-20.0000){\makebox(0,0)[lb]{$\begin{pmatrix}x_i \\ y_i \end{pmatrix}$}}%
\put(39.9000,-20.0000){\makebox(0,0)[lb]{$\begin{pmatrix}x_n \\ y_n \end{pmatrix}$}}%
\put(16.0000,-10.1000){\makebox(0,0)[lb]{$\begin{pmatrix}x'_1 \\ y'_1 \end{pmatrix}$}}%
\put(25.7000,-10.3000){\makebox(0,0)[lb]{$\begin{pmatrix}x'_i \\ y'_i \end{pmatrix}$}}%
\put(39.9000,-10.4200){\makebox(0,0)[lb]{$\begin{pmatrix}x'_n \\ y'_n \end{pmatrix}$}}%
\put(1.6000,-19.9000){\makebox(0,0)[lb]{$\begin{pmatrix}x_1 \\ y_1\end{pmatrix}$}}%
\put(33.8000,-24.5000){\makebox(0,0)[lb]{$K_i\begin{pmatrix}x_i \\ y_i \end{pmatrix}$}}%
%
\special{pn 20}%
\special{sh 1}%
\special{ar 1014 2002 10 10 0  6.28318530717959E+0000}%
\special{sh 1}%
\special{ar 2942 2002 10 10 0  6.28318530717959E+0000}%
\special{sh 1}%
\special{ar 4854 2002 10 10 0  6.28318530717959E+0000}%
\special{sh 1}%
\special{ar 4854 2002 10 10 0  6.28318530717959E+0000}%
%
\special{pn 8}%
\special{pa 2782 1442}%
\special{pa 2614 1530}%
\special{fp}%
\special{sh 1}%
\special{pa 2614 1530}%
\special{pa 2682 1517}%
\special{pa 2661 1505}%
\special{pa 2664 1481}%
\special{pa 2614 1530}%
\special{fp}%
\special{pa 3102 2570}%
\special{pa 3246 2482}%
\special{fp}%
\special{sh 1}%
\special{pa 3246 2482}%
\special{pa 3179 2500}%
\special{pa 3200 2510}%
\special{pa 3200 2534}%
\special{pa 3246 2482}%
\special{fp}%
\end{picture}%
  \end{center}
  \caption{(rational) tadpole junction}\label{fig1}
\end{figure}
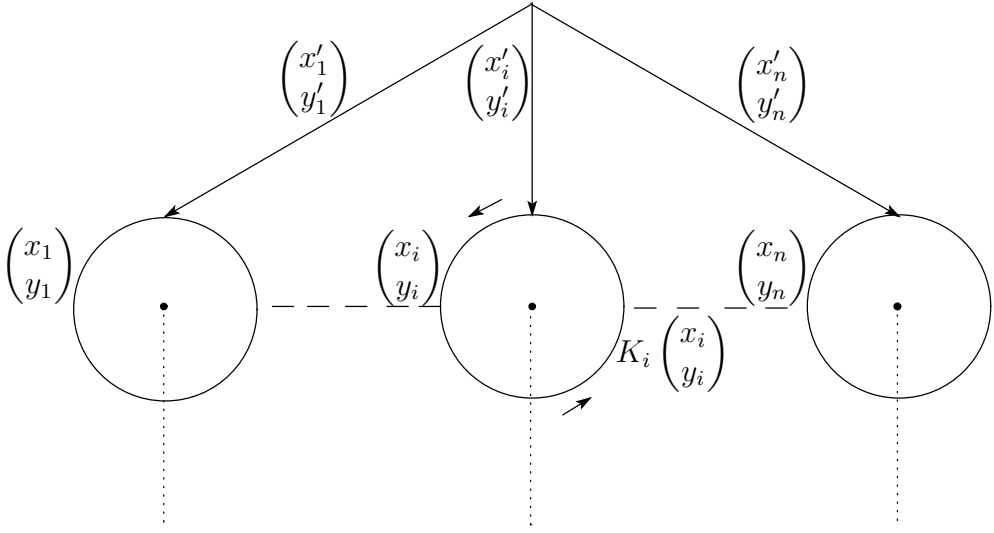

Consider $n$ open surfaces 
\begin{equation*}
 \rho'_1\left(\gamma_1\right),\rho'_2\left(\gamma_2\right),\dots,
 \rho'_n\left(\gamma_n\right) \quad
 \left(\gamma_i \in H_1\left(f^{-1}\left(u_0\right), 
 \bR \right)\right)
\end{equation*}
which have common boundary in $f^{-1}\left(u_0\right)$; glue them 
together to get a closed surface $J$ which is denoted by 
\begin{equation*}
 J = \rho'_1\left(\gamma_1\right)+\rho'_2\left(\gamma_2
 \right)+\dotsb+\rho'_n\left(\gamma_n\right).
\end{equation*}
We also use notations such as 
\begin{equation*}
\begin{split}
 J &= \left[x_1,y_1;x_2,y_2;\dots; 
      x_n,y_n\right] \\
   &= \left(x'_1,y'_1;x'_2,y'_2;\dots;x'_n,y'_n\right),
\end{split}
\end{equation*}
where 
 $\gamma_i = 
 x_i\alpha +y_i\beta \left(1 \le i \le n \right)$ and 
 $\gamma_i - {\rho'}_{i}^{*}\left(\gamma_i\right) = x'_i\alpha + y'_i\beta.$ 
Then we have 
\begin{equation}
\begin{pmatrix} x'_i \\ y'_i \end{pmatrix} =
\left(I - K_i\right)\begin{pmatrix} x_i \\ y_i \end{pmatrix} 
\label{sub1}
\end{equation}
\ake
For $J = \left[x_1,y_1;x_2,y_2;\dots;x_n,y_n\right] = 
\left(x'_1,y'_1;x'_2,y'_2;\dots;x'_n,y'_n\right)$, 
$J$ is said to be a \textit{tadpole junction} if all $x_i,y_i$ are integral 
and a \textit{rational tadpole junction} if all $x'_i,y'_i$ are integral 
(See Figure $\ref{fig1}$). 

Now we define the \textit{tadpole junction lattice} of $X$.
We put $T=\{ u_i \in S \ \vert \ f^{-1}(u_i) \ 
\textrm{is of type I} \}$ and 
$l = 2n - \vert T \vert$. 
If $u_i \in T$, we can assume that $f^{-1}\left(u_i\right)$ is of 
type $I_{n_i}$ and \begin{equation*}
 K_{i} = \begin{pmatrix} 
1+n_{i} p_{i} q_{i} & -n_{i} p_{i}^2 \\ 
n_{i} q_{i}^2 & 1+n_{i} p_{i} q_{i}
            \end{pmatrix}.
\end{equation*}
Then $x_i$ and $y_i$ are redundant as variables and we can use
variable $t_i$ defined by
\begin{equation*}
 t_{i} = 
 \begin{pmatrix} p_{i} \\ q_{i} \end{pmatrix} \times
 \begin{pmatrix} x_{i} \\ y_{i} \end{pmatrix},
\end{equation*}
where $\begin{pmatrix} a \\ b \end{pmatrix} \times 
\begin{pmatrix} c \\ d \end{pmatrix}$ means $ad-bc$.
It is clear that
\begin{equation*}
 \begin{pmatrix} x'_{i} \\ y'_{i} \end{pmatrix} = 
 n_{i}t_{i}\begin{pmatrix} p_{i} \\ q_{i} \end{pmatrix}. 
\end{equation*}
For $u_i \in T$ we will take $t_i$ as the independent variable instead of 
$x_i, y_i$ and we use the notation of $J$ such as
$J=[s_1;s_2;\ldots;s_n]$ where $s_i=(x_i,y_i)$ ($u_i \notin T$) and
$s_i=t_i$ ($u_i \in T$).

By the boundary condition, we have 
\begin{equation}
\sum_{u_i \notin T}
\left(I - K_i\right)\begin{pmatrix} x_i \\ y_i \end{pmatrix}+
\sum_{u_i \in T}n_{i}t_{i}
\begin{pmatrix} p_{i} \\ q_{i} \end{pmatrix} = 
\begin{pmatrix} 0 \\ 0 \end{pmatrix}.\label{pq1}
\end{equation}
Let
\begin{equation*}
W_{1} = \{J = \left[x_1,y_1;\dots;t_i;\dots;x_n,y_n \right] \in \bZ^{l} \mid 
x_i,y_i,t_i \>
\textrm{satisfies the condition} \>
\left(\ref{pq1}\right)\}
\end{equation*}
and 
\begin{equation*}
\langle \hspace{1.5mm}, \hspace{1.5mm}\rangle : W_{1} \times W_{1} 
\to \bZ ,
\quad\left(J_1, J_2\right) \mapsto \langle J_1, J_2\rangle 
\end{equation*}
be the symmetric bilinear pairing defined by the intersection number of 
the tadpole junctions $J_1$ and $J_2$. 
Then, for $J_1 = \left[x_1,y_1;\dots;t_i;\dots;x_n,y_n \right]$ and 
$J_2 = \left[u_1,v_1;\dots;w_i;\dots;u_n,v_n \right]$, we have 
\begin{equation}
\begin{split}
 \langle J_1, J_2\rangle 
 =&  \sum_{u_i \notin T} 
    \begin{pmatrix} x_i \\ y_i \end{pmatrix} 
    \times \left(I-K_i\right)\begin{pmatrix} u_i \\ v_i \end{pmatrix} 
 -\sum_{u_i \in T} n_{i} t_i w_i \\
 &+\sum_{1 \le i < j \le n}
 \left(I-K_i\right)\begin{pmatrix} x_i \\ y_i \end{pmatrix} \times 
 \left(I-K_j\right)\begin{pmatrix} u_j \\ v_j \end{pmatrix}.\label{seki1}
\end{split}
\end{equation}

Hence we can define the \textit{tadpole junction lattice} of $X$ by 
$\left(W_{1},\langle \hspace{1.5mm}, \hspace{1.5mm}\rangle\right)$. 
That is,

\begin{Pro}
{\normalfont\rmfamily
Let \,$\left(W_{1},\langle \hspace{1.5mm}, \hspace{1.5mm}\rangle_{K}\right)
(resp.\> ({\tilde{W}}_{1},\langle \hspace{1.5mm}, 
\hspace{1.5mm}\rangle_{\tilde{K}}))$ be a tadpole junction lattice for a 
global monodromy $K = \left(K_1,K_2,\dots,K_n\right)$
$(resp. \> \tilde{K} = ({\tilde{K}}_1,{\tilde{K}}_2,\dots,{\tilde{K}}_n))$. 
If $K \sim \tilde{K}$, then there exists an isomorphism between 
$\left(W_{1},\langle \hspace{1.5mm}, \hspace{1.5mm}\rangle_{K}\right)$ and 
$({\tilde{W}}_{1},\langle \hspace{1.5mm}, 
\hspace{1.5mm}\rangle_{\tilde{K}})$. 
}
\end{Pro}
\textit{Proof.} \, We can assume that $K$ transfer $\tilde{K}$ by the elementary 
transformation $R_i$. 
If $u_i \notin T$ and $u_{i+1} \notin T$, then we let the $4 \times 4$ matrix 
$P$
\begin{equation*}
 P = \begin{pmatrix}
     I_2 - K_i & I_2 \\
     K_{i+1}   & O 
     \end{pmatrix}.
\end{equation*}
If $u_i \in T$ and $u_{i+1} \notin T$, then we let the $3 \times 3$ matrix 
$P$
\begin{equation*}
 P = \begin{pmatrix}
     n_ip_i & 1 & 0 \\
     n_iq_i & 0 & 1 \\
     1                  & 0 & 0
     \end{pmatrix}.
\end{equation*}
If $u_i \notin T$ and $u_{i+1} \in T$, then we let the $3 \times 3$ matrix 
$P$
\begin{equation*}
 P = \begin{pmatrix}
     (-q_{i+1},p_{i+1})\left(I_2 -K_i\right) & 1 \\
     K_{i+1} & 0
     \end{pmatrix}.
\end{equation*}
If $u_i \in T$ and $u_{i+1} \in T$, then we let the $2 \times 2$ matrix 
$P$
\begin{equation*}
 P = \begin{pmatrix}
     {n_i}(p_{i+1}q_{i}-p_{i}q_{i+1}) & 1 \\
     1 & 0
     \end{pmatrix}.
\end{equation*}
Let the $l \times l$ matrix $P_i$ 
\begin{equation*}
 P_i = \begin{pmatrix} 
      \> 1 &        &   &   &   &        &   \\
         & \ddots &   &   &   &        &   \\
         &        & 1 &   &   &        &   \\
         &        &   & P &   &        &   \\
         &        &   &   & 1 &        &   \\
         &        &   &   &   & \ddots &   \\
         &        &   &   &   &        & 1\> \\
       \end{pmatrix}.
\end{equation*}
It is clear that 
\begin{equation*}
 det(P_i) = \pm 1.
\end{equation*}
For $J = \left[x_1,y_1;\dots;x_n,y_n\right] \in W_1$, we let 
$\tilde{J} = {}^{t}P_i J$. 
Then we have 
$\tilde{J} \in W_2$ and 
\begin{equation*}
\langle J, J \rangle_{K} = 
\langle \tilde{J}, \tilde{J} \rangle_{\tilde{K}}.
\end{equation*}
\textit{Q.E.D.}
\ake

Next we define the \textit{rational tadpole junction lattice} of $X$.
Consider a tadpole junction $J' = \left(x'_1,y'_1;\dots;x'_n,
y'_n \right) \in \bZ^{l}$.\\
If $u_i \in T$, we let 
\begin{equation}
t'_i = n_i t_i. \label{sub2}
\end{equation}
By substituting $\left(\ref{sub1}\right)$ and $\left(\ref{sub2}\right)$ 
to $\left(\ref{pq1}\right)$, 
we have 
\begin{equation}
\sum_{u_i \notin T}
\begin{pmatrix} x'_i \\ y'_i \end{pmatrix}+
\sum_{u_i \in T}t'_{i}
\begin{pmatrix} p_{i} \\ q_{i} \end{pmatrix} = 
\begin{pmatrix} 0 \\ 0 \end{pmatrix}. \label{pq2}
\end{equation}
Let 
\begin{equation*}
W_{2} = \{J' = \left( x'_1,y'_1;\dots;t'_i;\dots;x'_n,y'_n \right) \in 
\bZ^{l} \mid 
x'_i,y'_i,t'_i \>
\textrm{satisfies the condition} \>
\left(\ref{pq2}\right)\}
\end{equation*}
and 
\begin{equation*}
\langle \hspace{1.5mm}, \hspace{1.5mm}\rangle : W_{2} \times W_{2} \to \bQ ,
\quad\left(J'_1, J'_2\right) \mapsto \langle J'_1, J'_2\rangle
\end{equation*}
be the symmetric bilinear pairing defined by the intersection number of 
the rational tadpole junctions $J'_1$ and $J'_2$. 
Then we have 
\begin{equation}
\begin{split}
 \langle J'_1, J'_2\rangle 
 =&  \sum_{u_i \notin T} 
    \left(I-K_i\right)^{-1} 
    \begin{pmatrix} x'_i \\ y'_i \end{pmatrix} 
    \times \begin{pmatrix} u'_i \\ v'_i \end{pmatrix} 
 -\sum_{u_i \in T}\frac{t'_i w'_i}{n_{i}} \\
 &+\sum_{1 \le i < j \le n}
 \begin{pmatrix} x'_i \\ y'_i \end{pmatrix} \times 
 \begin{pmatrix} u'_j \\ v'_j \end{pmatrix}  \label{seki2}
\end{split}
\end{equation}
for $J'_1 = \left( x'_1,y'_1;\dots;t'_i;\dots;x'_n,y'_n \right)$ and 
$J'_2 = \left( u'_1,v'_1;\dots;w'_i;\dots;u'_n,v'_n \right)$.
We can define 
the \textit{rational tadpole junction lattice} of $X$ by 
$\left(W_{2},\langle \hspace{1.5mm}, \hspace{1.5mm}\rangle\right)$. 
That is, 

\begin{Pro}
{\normalfont\rmfamily
Let \,$\left(W_{2},\langle \hspace{1.5mm}, \hspace{1.5mm}\rangle_{K}\right)
(resp.\> ({\tilde{W}}_{2},\langle \hspace{1.5mm}, 
\hspace{1.5mm}\rangle_{\tilde{K}}))$ be a tadpole junction lattice for a 
global monodromy $K = \left(K_1,K_2,\dots,K_n\right)$
$(resp. \> \tilde{K} = ({\tilde{K}}_1,{\tilde{K}}_2,\dots,{\tilde{K}}_n))$. 
If $K \sim \tilde{K}$, then there exists an isomorphism between 
$\left(W_{2},\langle \hspace{1.5mm}, \hspace{1.5mm}\rangle_{K}\right)$ and 
$({\tilde{W}}_{2},\langle \hspace{1.5mm}, 
\hspace{1.5mm}\rangle_{\tilde{K}})$. 
}
\end{Pro}
\textit{Proof.} \, We can assume that $K$ transfer $\tilde{K}$ by the elementary 
transformation $R_i$. 
If $u_i \notin T$ and $u_{i+1} \notin T$, then we let the $4 \times 4$ matrix 
$Q$
\begin{equation*}
 Q = \begin{pmatrix}
     I_2 - K_{i+1} & I_2 \\
     K_{i+1}   & O 
     \end{pmatrix}.
\end{equation*}
If $u_i \in T$ and $u_{i+1} \notin T$, then we let the $3 \times 3$ matrix 
$Q$
\begin{equation*}
 Q = \begin{pmatrix}
     \left(I_2-K_{i+1}\right) 
     \begin{pmatrix}p_i \\ q_i \end{pmatrix} & I_2 \\
     1                  & O_{1,2}
     \end{pmatrix}.
\end{equation*}
If $u_i \notin T$ and $u_{i+1} \in T$, then we let the $3 \times 3$ matrix 
$Q$
\begin{equation*}
 Q = \begin{pmatrix}
     -n_{i+1}q_{i+1},\> n_{i+1}p_{i+1}  & 1 \\
     K_{i+1} & \begin{matrix}0 \\ 0\end{matrix}
     \end{pmatrix}.
\end{equation*}
If $u_i \in T$ and $u_{i+1} \in T$, then we let the $2 \times 2$ matrix 
$Q$
\begin{equation*}
 Q = \begin{pmatrix}
     {n_{i+1}}(p_{i+1}\bar{y_{i}}-p_{i}q_{i+1}) & 1 \\
     1 & 0
     \end{pmatrix}.
\end{equation*}
Let the $l \times l$ matrix $Q_i$ 
\begin{equation*}
 Q_i = \begin{pmatrix} 
      \> 1 &        &   &   &   &        &   \\
         & \ddots &   &   &   &        &   \\
         &        & 1 &   &   &        &   \\
         &        &   & Q &   &        &   \\
         &        &   &   & 1 &        &   \\
         &        &   &   &   & \ddots &   \\
         &        &   &   &   &        & 1\> \\
       \end{pmatrix}.
\end{equation*}
It is clear that 
\begin{equation*}
 det(Q_i) = \pm 1.
\end{equation*}
For $J'= \left(x'_1,y'_1;\dots;x'_n,y'_n\right) \in W_2$, we let 
$\tilde{J'} = {}^{t}Q_i J'$. 
Then we have 
$\tilde{J'} \in W_2$ and 
\begin{equation*}
\langle J', J' \rangle_{K} = 
\langle \tilde{J'}, \tilde{J'} \rangle_{\tilde{K}}.
\end{equation*}
\textit{Q.E.D.}
\ake

\subsection{Null junction and rational null junction}\label{null}
Let $\rho_0$ be an path which is homotopic to constant mapping $id_{u_0}$.
We define a \textit{null junction} to be a tadpole junction which is 
homologous to  
$\rho_0\left(\gamma\right)$ 
$\left(\gamma' \in H_1\left(f^{-1}(u_0), \bZ\right)\right)$ 
and let 
\begin{equation*}
 U_{1} = \{J = \left[x_1,y_1;\dots;x_n,y_n \right] 
             = \left(x'_1,y'_1;\dots;x'_n,y'_n \right) 
             \in \bZ^{l} \mid 
           J \>\> \textrm{is a null junction}\},
\end{equation*}
which is called the \textit{null junction \bZ-module}. 

   \begin{Lem}\label{base1}
{\normalfont\rmfamily
The null junction \bZ-module $U_1$ is generated by 
$J_a=\left[x_1,y_1;\dots;x_n,y_n \right] \in \bZ^{2n}$ and 
$J_b=\left[u_1,v_1;\dots;u_n,v_n \right] \in \bZ^{2n}$, where 
\begin{equation*}
 \begin{pmatrix} x_i \\ y_i \end{pmatrix} = 
 K_{i-1}K_{i-2}\cdots K_1 \begin{pmatrix} 1 \\ 0 \end{pmatrix}, \>\>
 \begin{pmatrix} u_i \\ v_i \end{pmatrix} = 
 K_{i-1}K_{i-2}\cdots K_1 \begin{pmatrix} 0 \\ 1 \end{pmatrix} \>\>\>\>
 \left(1 \le i \le n \right).
\end{equation*}
}
\end{Lem}
\textit{Proof.} \, Because $\rho_0 = \rho_1 \rho_2 \cdots \rho_n$, by using 
$\left(\ref{rho}\right)$ repeatedly, we obtain 
\begin{equation*}
\rho_0\left(\gamma\right) = \rho_1\left(\gamma\right) + \cdots + 
\rho_i\left(\rho_{i-1}^{*}\cdots \rho_{1}^{*} (\gamma)\right) +\cdots + 
\rho_n\left(\rho_{n-1}^{*}\cdots \rho_{1}^{*} (\gamma)\right).
\end{equation*}
We put $\gamma = p\alpha + q\beta \left(p,q \in \bZ \right)$. 
We have $\rho_0\left(\gamma\right) = \rho_0\left(\alpha\right) + 
\rho_0\left(\beta\right)$. 
And, by the definition of the tadpole junction, we have 
\begin{gather*}
J_a = \rho_1\left(\alpha\right) + \cdots + 
\rho_i\left(\rho_{i-1}^{*}\cdots \rho_{1}^{*} (\alpha)\right) +\cdots + 
\rho_n\left(\rho_{n-1}^{*}\cdots \rho_{1}^{*} (\alpha)\right) \\
J_b = \rho_1\left(\beta\right) + \cdots + 
\rho_i\left(\rho_{i-1}^{*}\cdots \rho_{1}^{*} (\beta)\right) +\cdots + 
\rho_n\left(\rho_{n-1}^{*}\cdots \rho_{1}^{*} (\beta)\right).
\end{gather*}
We have thus proved the lemma.
\textit{Q.E.D.}
\ake

Similarly we define a \textit{rational null junction} to be a rational 
tadpole junction 
which is homologous to $\rho_0\left(\gamma'\right)$ 
$\left(\gamma' \in H_1\left(f^{-1}(u_0), \bQ\right)\right)$ 
and let 
\begin{equation*}
 U_{2} = \{J' = \left(x'_1,y'_1;\dots;x'_n,y'_n \right) \in \bZ^{l} \mid 
           J' \>\> \textrm{is a rational null junction}\},
\end{equation*}
which is called the \textit{rational null junction \bZ-module}.\\
For $i\left(1\le i \le n\right)$, if $u_i \in T$, then let 
$N_i = n_i\begin{pmatrix} -q_i , p_i \end{pmatrix}$ else 
$N_i = I - K_i$.
Let the $l \times 2$ matrix N 
\begin{equation*}
 N = \begin{pmatrix}
     N_1\\
     \vdots\\
     N_n
     \end{pmatrix}, 
\end{equation*} 
and the elementary divisor type of $N$ be $\left(n_1,n_2\right)$. 
That is, there are a $l \times l$ unimodular matrix $A$ and 
a $2 \times 2$ unimodular matrix $B = \left(b_1,b_2\right)$ such that 
\begin{equation}
ANB = \begin{pmatrix} 
       n_1 & 0 \\
       0   & n_2\\
       0   & 0\\
       \vdots & \vdots\\
       0 & 0
       \end{pmatrix}.\label{GMH}
\end{equation}
\begin{Lem}\label{base2}
{\normalfont\rmfamily
The rational null junction \bZ-module $U_2$ is generated by 
$J'_a=\left(x'_1,y'_1;\dots;x'_n,y'_n \right)$ and 
$J'_b=\left(u'_1,v'_1;\dots;u'_n,v'_n \right)$, where 
\begin{equation}
 \begin{pmatrix} x'_i \\ y'_i \end{pmatrix} = 
 \frac{1}{n_1}\left(I-K_i\right)K_{i-1}\cdots K_1 b_1, \>\>
 \begin{pmatrix} u'_i \\ v'_i \end{pmatrix} = 
 \frac{1}{n_2}\left(I-K_i\right)K_{i-1}\cdots K_1 b_2\>\>\>\>
 \left(1 \le i \le n \right). \label{base2siki}
\end{equation}
}
\end{Lem}
\textit{Proof.} \, Let $r J_a + s J_b \left(r,s \in \bQ\right)$ be a 
rational null junction, and put $\>{}^{t}\left(r',s'\right) = B^{-1}
\>\>\>{}^{t}\left(r,s\right)$. 
For $i\left(1\le i \le n\right)$, if $u_i \in T$, then let 
\begin{equation*}
 M_i = n_i\left(
       \begin{pmatrix} p_i \\ q_i \end{pmatrix} \times 
       K_{i-1} \cdots K_1 \begin{pmatrix} 1 \\ 0 \end{pmatrix},
       \begin{pmatrix} p_i \\ q_i \end{pmatrix} \times 
       K_{i-1} \cdots K_1 \begin{pmatrix} 0 \\ 1 \end{pmatrix}
                      \right)
\end{equation*}
else 
\begin{equation*}
 M_i = \left(I - K_i\right)K_{i-1} \cdots K_1.
\end{equation*}
Let the $l \times 2$ matrix $M$
\begin{equation*}
 M = \begin{pmatrix}
     M_1\\
     \vdots\\
     M_n
     \end{pmatrix}.
\end{equation*} 
By the definition of the rational tadpole junction, we get 
$M\>{}^{t}\left(r,s\right) \in \bZ^{l}$. 
This is equivalent to 
the condition $N\>{}^{t}\left(r,s\right) \in \bZ^{l}$, and 
which is equivalent to the condition 
$ANB\>{}^{t}\left(r',s'\right) \in \bZ^{l}$ because $A$ is unimodular. 
By $\left(\ref{GMH}\right)$ it means $n_1 x', n_2 y' \in \bZ$, 
and which means $n_1 x, n_2 y \in \bZ$ because $B$ is unimodular. 
That is, $x = a / n_1, y = b / n_2 \left(a,b \in \bZ\right)$. 
Therefore $J'_a = \left(1/n_1\right) J_a, J'_b = \left(1/n_2\right) J_b$. 
By Lemma \ref{base1} we have $\left(\ref{base2siki}\right)$. 
\textit{Q.E.D.}
\ake
{}From the definition, we can identify the tadpole junction with the 
element of $H_2\left(X, \bZ\right)$ by the class mapping 
\begin{equation*}
 h : W_1 \rightarrow V^{\perp} \subset L_k,\>\> J \mapsto [J]. 
\end{equation*}
\begin{Lem}\label{Kerh}
{\normalfont\rmfamily 
\begin{equation*}
 \textrm{Ker} \> h = U_1 
\end{equation*}
}
\end{Lem}
\textit{Proof.} \, It is enough to show $\textrm{Ker} \> h \subset U_1$. 
Let $J = \left[x_1,y_1;\dots;x_n,y_n\right] \in \bZ^{2n}$ be the element of 
$\textrm{Ker} \> h$. 
For $i\left(1\le i \le n-1\right)$, define the tadpole junction 
$\check{J}_i,\hat{J}_i$ by 
\begin{gather*}
\check{J}_i = [\overbrace{0,0}^{1};\overbrace{0,0}^{2};\dots;
\overbrace{1,0}^{i};\overbrace{-1,0}^{i+1};\dots;
\overbrace{0,0}^{n-1};\overbrace{0,0}^{n}] \\ 
\hat{J}_i = [\overbrace{0,0}^{1};\overbrace{0,0}^{2};\dots;
\overbrace{0,1}^{i};\overbrace{0,-1}^{i+1};\dots;
\overbrace{0,0}^{n-1};\overbrace{0,0}^{n}].
\end{gather*} 
We have $\langle J,\check{J}_i \rangle = 0$ and $\langle J,\hat{J}_i 
\rangle = 0$ because $J$ is the element of $Ker \> h$. 
On the other hand, by an easy calculation, we obtain 
\begin{equation*}
\langle J,\check{J}_i \rangle = 
\left(
\begin{pmatrix}x_i \\ y_i \end{pmatrix} - 
K_{i-1}\begin{pmatrix}x_{i-1} \\ y_{i-1}\end{pmatrix}
\right) \times \begin{pmatrix}-1 \\ 0\end{pmatrix}, \>\> 
\langle J,\hat{J}_i \rangle = 
\left(
\begin{pmatrix}x_i \\ y_i \end{pmatrix} - 
K_{i-1}\begin{pmatrix}x_{i-1} \\ y_{i-1}\end{pmatrix}
\right) \times \begin{pmatrix}0 \\ -1\end{pmatrix}.
\end{equation*} 
Therefore we have 
$\begin{pmatrix}x_i \\ y_i \end{pmatrix} = 
K_{i-1}\begin{pmatrix}x_{i-1} \\ y_{i-1}\end{pmatrix}$. 
And then we have $\begin{pmatrix} x_i \\ y_i \end{pmatrix} = 
K_{i-1}\cdots K_1 \begin{pmatrix} x_1 \\ y_1 \end{pmatrix}$. 
\vskip 1mm \noindent
By Lemma \ref{base1} we have $J = x_1 J_a + y_1 J_b \in U_1$. 
\textit{Q.E.D.}
\ake

\subsection{Brane configuration and junction}
In this section we define the brane configuration and junction lattice and 
prove the necessary proposition and the lemmas.
\ake
It is well known that 
the local monodromy matrix $K$ of the singular fibre $F$ of type $I_{1}$ 
\begin{equation*}
 K = \begin{pmatrix} 1+pq & -p^{2} \\ q^{2} & 1-pq
     \end{pmatrix}
\end{equation*}
if the vanishing cycle of $F$ is $p \alpha + q \beta$. 
Then we denote the singular fibre $F$ by $\bX_{\left[p,q\right]}$ and 
the matrix $K$ by $K_{\left[p,q\right]}$. 
In particular we denote $\bX_{\left[1,0\right]}$, 
$\bX_{\left[1,-1\right]}$, $\bX_{\left[1,1\right]}$ by $\bA$, $\bB$, $\bC$ 
and $K_{\left[1,0\right]}$, $K_{\left[1,-1\right]}$, 
$K_{\left[1,1\right]}$ by $K_A$, $K_B$, $K_C$ respectively. 
We call the finite sequence of $\bX_{\left[p,q\right]}$ the
\textit{brane configuration}. 
Let $X = \bX_{\left[p_1,q_1\right]} \cdots 
\bX_{\left[p_n,q_n\right]}$. By $\left(\ref{elem}\right)$, we can define 
the \textit{elementary transformations} 
$T_i,\, T_i^{-1}$ of brane configurations
\begin{equation}
\begin{split}
& T_i\left(
\bX_{\left[p_1,q_1\right]}\cdots
\bX_{\left[p_i,q_i\right]}\bX_{\left[p_{i+1},
q_{i+1}\right]}
\cdots\bX_{\left[p_{n},q_{n}\right]}
\right) = 
\bX_{\left[p_1,q_1\right]}\cdots
\bX_{\left[p_{i+1},q_{i+1}\right]}\bX_{\left[r,s\right]}
\cdots\bX_{\left[p_{n},q_{n}\right]} \\
& T_{i}^{-1}\left(
\bX_{\left[p_1,q_1\right]}\cdots
\bX_{\left[p_i,q_i\right]}\bX_{\left[p_{i+1},
q_{i+1}\right]}
\cdots\bX_{\left[p_{n},q_{n}\right]}
\right) = 
\bX_{\left[p_1,q_1\right]}\cdots
\bX_{\left[t,u\right]}\bX_{\left[p_i,q_i\right]}
\cdots\bX_{\left[p_{n},q_{n}\right]}, \label{T+-}
\end{split}
\end{equation}
where 
\begin{gather*}
\left[r,s\right] = \left[p_i,q_i\right] + 
\left(p_i q_{i+1} - p_{i+1} q_{i}\right)
\left[p_{i+1},q_{i+1}\right] \\
\left[t,u\right] = \left[p_{i+1},q_{i+1}\right] + 
\left(p_i q_{i+1} - p_{i+1} q_{i}\right)
\left[p_{i},q_{i}\right].
\end{gather*}
If there exists a finite sequence of elementary transformation which 
transfer $\bX_{\left[p_1,q_1\right]} \cdots 
\bX_{\left[p_n,q_n\right]}$ 
to $\bX_{\left[r_1,s_1\right]} \cdots 
\bX_{\left[r_n,s_n\right]}$, then 
we identify 
$\bX_{\left[p_1,q_1\right]} \cdots 
\bX_{\left[p_n,q_n\right]}$ with 
$\bX_{\left[r_1,s_1\right]} \cdots 
\bX_{\left[r_n,s_n\right]}$ and 
denote $\bX_{\left[p_1,q_1\right]} \cdots 
\bX_{\left[p_n,q_n\right]}\sim 
\bX_{\left[r_1,s_1\right]} \cdots 
\bX_{\left[r_n,s_n\right]}$.\\
Now we define the \textit{junction lattice} of the brane configuration $X$.
We define the $\bZ$-module $Z_{X}$ of rank $n-2$ by
\begin{equation*}
 Z_{X} = \left\{Q = \left(Q_1,\dots,Q_n\right)\in \bZ^{n} \> \bigg| \> 
       Q_1\begin{pmatrix}p_1 \\ q_1\end{pmatrix} + \cdots + 
       Q_n\begin{pmatrix}p_n \\ q_n\end{pmatrix} = 
       \begin{pmatrix}0 \\ 0\end{pmatrix}
     \right\},
\end{equation*}
and the symmetric bilinear $\left(Q, Q\right)$ pairing by 
\begin{equation*}
 \left(Q, Q\right) = -\sum_{i=1}^{n}Q_{i}^{2} + 
 \sum_{i < j}\left(p_i q_j - p_j q_i\right)Q_i Q_j.
\end{equation*}
We call the lattice 
$\left(Z_{X},\left(\hspace{1.5mm}, \hspace{1.5mm}\right) \right)$ 
\textit{junction lattice} of $X$ and call the element of 
$Z_X$ \textit{junction}. 
This definition is well-defined \cite{FYY}. 
We remark that the junction $Q = \left(Q_1,\dots,Q_n\right)$ can be 
geometrically interpreted as the tadpole junction 
$J = \left(Q_1(p_1,q_1);\dots;Q_n(p_n,q_n)\right)$. 
By an easy calculation we obtain 
\begin{equation}
\left(Q,Q\right) = \langle J,J \rangle. \label{koutensuu}
\end{equation}
\textbf{Theorem 4.}(\cite{M})\\
Let $\bX_{[p_{1},q_{1}]}\cdots
\bX_{[p_{l},q_{l}]}$ be the brane configuration which 
satisfies 
\begin{equation}
K_{[p_{l},q_{l}]} \cdots K_{[p_{1},q_{1}]} = I. 
\end{equation} 
Then there exists a natural number $k$ such that $l = 12k$ and 
\begin{equation}
 \bX_{[p_{1},q_{1}]}\cdots\bX_{[p_{12k},q_{12k}]} \sim 
 \underbrace{\bA^{8}\bB\bC\bB\bC}_{12}
 \cdots \underbrace{\bA^{8}\bB\bC\bB\bC}_{12} = 
 \left\{\bA^{8}\bB\bC\bB\bC\right\}^{k}.
\label{A8BCBC}
\end{equation}
\ake
\textbf{Theorem 5.}(\cite{DHIZ})\\
Let $X$ be the brane configuration $\bX_{[p_{1},q_{1}]}\cdots 
\bX_{[p_{n},q_{n}]}$ and $K_X$ be the matrix 
$K_{[p_{1},q_{1}]}\cdots K_{[p_{n},q_{n}]}$. 
Let $M_X$ be the Gram matrix of the junction lattice $Z_X$ of $X$ and $d$ 
is the G.C.D of $x_i y_j - x_j y_i \left(1 \le i < j \le n\right)$.
Except for the case $X = \bX_{\left[p,q\right]}^{n}$, we have 
\begin{equation*}
\textrm{det} \left(M_X\right) = \frac{\textrm{det} 
\left(I - K_X\right)}{d^{2}}\>. 
\end{equation*}
\ake

\begin{Pro}\label{jlattice}
{\normalfont\rmfamily
Let $X = \bX_{[p_{1},q_{1}]}\cdots
\bX_{[p_{12k},q_{12k}]}$ be the brane configuration which 
satisfies the condition $\left(\ref{A8BCBC}\right)$ 
and $Z_{X}$ the junction lattice of $X$.
Let $U$ be a $\bZ$-module generated by 
$Q_a = \left(a_1,\dots,a_{12k}\right)$ and 
$Q_b = \left(b_1,\dots,b_{12k}\right)$, where 
\begin{equation*}
a_i = \begin{pmatrix}p_i \\ q_i\end{pmatrix} \times 
K_{\left[p_{i-1},q_{i-1}\right]}\cdots 
K_{\left[p_{1},q_{1}\right]} 
\begin{pmatrix}1 \\ 0\end{pmatrix}, \>\>
b_i = \begin{pmatrix}p_i \\ q_i\end{pmatrix} \times 
K_{\left[p_{i-1},q_{i-1}\right]}\cdots 
K_{\left[p_{1},q_{1}\right]} 
\begin{pmatrix}0 \\ 1\end{pmatrix}. 
\end{equation*}
Then the quotient lattice $Z_{X}/U$ is isomorphic to $L_k$. 
}
\end{Pro}
\textit{Proof.} \, \\
\textit{Claim} $1.$ \, $(Q,Q)$ is even for arbitrary $Q = 
\left(Q_1,\dots,Q_{12k}\right) \in Z_{X}$. 
\ake
\textit{proof.} \, 
By induction, we can prove the following lemma. 
\begin{Lem}
{\normalfont\rmfamily
Let  $Q = \left(Q_1,\dots,Q_{12k}\right)$ is an element of $\bZ^{12k}$. 
If 
\begin{equation*}
Q_1 \begin{pmatrix} p_1 \\ q_1 \end{pmatrix} + \cdots + 
Q_{12k} \begin{pmatrix} p_{12k} \\ q_{12k} \end{pmatrix} \equiv 
\begin{pmatrix} 0 \\ 0 \end{pmatrix} \quad \left( mod \>\> 2 \right),
\end{equation*}
then $\left(Q,Q\right)$ is even, else $\left(Q,Q\right)$ is odd. 
}
\end{Lem}
Claim $1.$ follows from this lemma and the definition of $Z_{X}$. 
\textit{q.e.d.}
\ake
For a brane configuration $Y$, we denote the positive signature, the negative 
signature and the nullity of the Gram matrix of $Z_{Y}$ by 
$\sigma_{+}\left(Y\right)$, $\sigma_{-}\left(Y\right)$ and 
$\sigma_{0}\left(Y\right)$. 
For a matrix $A = \begin{pmatrix}
                  a_{1,1} & \cdots & a_{1,n} \\
                  \vdots  & \ddots & \vdots  \\
                  a_{n,1} & \cdots & a_{n,n}
                  \end{pmatrix}$, 
we denote $\begin{pmatrix}
               a_{1,1} & \cdots & a_{1,i} \\
               \vdots  & \ddots & \vdots  \\
               a_{i,1} & \cdots & a_{i,i}
               \end{pmatrix}$ by $A_i$ $\left(1 \le i \le n \right)$. \\
\ake
\textit{Claim} $2.$ \, $\sigma_{+}\left(X\right) = 2k-2,\>\> 
\sigma_{-}\left(X\right) = 10k-2,\>\> \sigma_{0}\left(X\right) = 2.$
\ake
\textit{proof.} \, 
By $\left(\ref{T+-}\right)$, we have 
\begin{equation*}
X \sim X_1 X_2 
\end{equation*}
where $X_1 = \bA^{8k}\bB^{2k}$ and 
$X_2 = \bX_{[4k-3,1-4k]}\bX_{[4k-5,3-4k]}\cdots\bX_{[1,-3]}\bC$. \\
First we show that $\sigma_{-}\left(X\right) \ge 10k-2.$ \\
In fact we have $Z_{X_1} \cong \left(- A_{8k-1}\right) 
\oplus \left(- A_{2k-1}\right)$, and then $\sigma_{-}\left(X\right) \ge 
\sigma_{-}\left(X_1\right) = 10k-2.$ \\
Next we show that $\sigma_{+}\left(X\right) \ge 2k-2.$ \\
It is enough to show that $Z_{X_2}$ is positive definite. 
Let $Y_n = \bX_{[2n-3,1-2n]}\bX_{[2n-5,3-2n]}\cdots\bX_{[1,-3]}\bC$ and 
$K_{Y_n} := \begin{pmatrix} a_n & b_n \\
                            c_n & d_n 
            \end{pmatrix}
         = K_{[-1,-1]}\cdots K_{[2n-5,3-2n]}K_{[2n-3,1-2n]}$ 
$\left(n \ge 3\right)$. 
Then by induction we have 
\begin{equation}
\begin{split}
a_n = (-1)^{n+1}(4n^2-n-1),\>\> b_n = (-1)^{n+1}(4n^2-5n),\\
c_n = (-1)^{n+1}n,\>\> d_n = (-1)^{n+1}(n-1). \label{abcd}
\end{split}
\end{equation}
Let $b_{n,i} = \underbrace{(0,\dots,0,\overbrace{1,0,\dots,0,1-i,2-i}^i)}
_n$ $\left(3 \le i \le n\right)$. It is easy to see that 
$B_i = \langle b_{i,3},\dots,b_{i,i} \rangle$ is the basis of $Z_{Y_i}$. 
Let $M_{Y_i}$ is the Gram matrix of $Z_{Y_i}$ about the basis $B_i$ 
$\left(3 \le i \le n\right)$. Then by Theorem $5.$ and 
$\left(\ref{abcd}\right)$ we have 
\begin{equation*}
\begin{split}
det\left(\left(M_{Y_n}\right)_i\right)&= det \left(M_{Y_i}\right) \\
                                      &= \left(4i^2-2+(-1)^{i}2\right)/16
 > 0 \>\> \left(3 \le i \le n\right).
\end{split}
\end{equation*}
Therefore it follows that $Z_{Y_n}$ is positive definite. 
We see that $Z_{X_2}$ is positive definite because $X_2 = Y_{2k}$. \\
Lastly we have $\sigma_{0}\left(X\right) \ge 2$ by Lemma 
$\ref{base1}$ and Lemma $\ref{Kerh}$. 
Therefore we have proved Claim $2.$. \textit{q.e.d.}
\ake
By $\left(\ref{T+-}\right)$, we have 
\begin{equation*}
X \sim X_3 \bX_{[3,1]}\bA ,
\end{equation*}
where $X_3 = \left(\bA^8\bB\bC\bB\bC\right)^{n-1}\bA^7\bB\bC^2$. 
Let $M_{X_3}$ be the Gram matrix of $Z_{X_3}$ 
and $C = \langle c_1,\dots,c_{12k-4}\rangle$ the basis of $Z_{X_3}$. 
By Theorem $5.$ we have $det\left(M_{X_3}\right) = 1$. 
Hence it follows from Lemma $\ref{base1}$ and Lemma $\ref{Kerh}$ 
that $\langle c_1,\dots,c_{12k-4},Q_a,Q_b\rangle$ is the basis of $Z_{X}$. 
Therefore $Z_{X}/U$ is isomorphic to $Z_{X_3}$. 
By Claim $1.$ and Claim $2.$ we see that $Z_{X_3}$ is even unimodular 
lattice and the signature of $Z_{X_3}$ is $\left(2k-2,10k-2\right)$. 
We have $Z_{X_3} \cong L_k$ by lattice theory \cite{Se}. 
We have thus proved the proposition. 
\textit{Q.E.D.}
\ake

In the following we describe the confluence of singular fibres by the 
brane configuration. 
Consider the case that the singular fibre $F$ is of type $I_{0}^{*}$ for 
instance. 
In this case it is known that the local monodromy matrix of 
$F$ is $-I$ \cite{Kod}. 
The brane configuration $\bX_{F} = \bA^{4}\bB\bC$ corresponds to $F$ 
because $-I = K_{\left[1,1\right]}K_{\left[1,-1\right]}
\left(K_{\left[1,0\right]}\right)^{4}$. 
We remark that this kind of expression of the monodoromy matrix 
as a confluence is not unique. However, all such expressions are
equivalent (see for insetance \cite{Na}.)
And then it turns 
out that the junction lattice $Z_{X_{F}}$ is $Z_{X_{F}} = -D_4$ \cite{DZ}, 
which corresponds to 
the fact that $-T_{v}$ of the singular fibre of type $I_{0}^{*}$ is $D_4$. 
In general, the brane configuration $X_{F}$ and the junction lattice 
$Z_{X_F}$ of the singular fibre $F$ can be summarized in Table 
$\ref{table2}$ \cite{DZ}. 
\begin{table}
\begin{center}
\begin{tabular}{|c|c|c|} \hline
\multicolumn{1}{|c|}{the type of $F$ } & brane configuration $X_{F}$ 
& junction lattice $Z_{X_{F}}$ \\ \hline
$I_n\left(1\le n\right)$ & $\bA^{n}$ & $-A_{n-1}\left(A_0 = 0\right)$ \\ \hline
\setcounter{part}{2} $\thepart$ & $\bA\bC$ & $0$ \\ \hline
\setcounter{part}{3} $\thepart$ & $\bA^{2}\bC$ & $-A_1$ \\ \hline
\setcounter{part}{4} $\thepart$ & $\bA^{3}\bC$ &$-A_2$ \\ \hline
$I_{n}^{*}\left(0\le n\right)$ & $\bA^{n+4}\bB\bC$ & $-D_{n+4}$ \\ \hline
\setcounter{part}{2} ${\thepart}^{*}$ & $\bA^{7}\bB\bC^{2}$ & $-E_8$ \\ \hline
\setcounter{part}{3} ${\thepart}^{*}$ & $\bA^{6}\bB\bC^{2}$ & $-E_7$ \\ \hline
\setcounter{part}{4} ${\thepart}^{*}$ & $\bA^{5}\bB\bC^{2}$ & $-E_6$ \\ \hline
\end{tabular}
\end{center}
\caption{the brane configuration $X_{F}$ and the junction lattice $Z_{X_F}$ of the singular fibre $F$}
\label{table2}
\end{table}
\ake
\ake
\indent Consider the elliptic surface $X_k$ such that the global monodromy is 
$K = \left(K_1,\dots,K_n\right)$. 
Let the brane configuration $X_{F_{u_i}}$ which correspond to the singular 
fibre $F_{u_i}$ 
\begin{equation*}
X_{F_{u_i}} = \bX_{\left[p_{i,1},q_{i,1}\right]}\cdots \bX_{
[p_{i,m_i},q_{i,m_i}]} \quad \left(1 \le i \le n\right),
\end{equation*} 
and then we have 
\begin{equation*}
K_i = K_{[p_{i,m_i},q_{i,m_i}]}\cdots 
K_{\left[p_{i,1},q_{i,1}\right]}.
\end{equation*}
We associate the brane configuration 
$X = \left(X_{F_{u_1}}\right) \cdots \left(X_{F_{u_n}}\right)$ with 
the global monodromy $K$, and then we call $X$ the 
\textit{brane configuration} of $X_k$. 
We put $\tilde{X} = X_{F_{u_1}} \cdots X_{F_{u_n}} = 
\bX_{\left[p_{1,1},q_{1,1}\right]}\cdots \bX_{
[p_{1,m_1},q_{1,m_1}]} \cdots \cdots \cdots
\bX_{\left[p_{n,1},q_{n,1}\right]}\cdots \bX_{
[p_{n,m_n},q_{n,m_n}]}$. 
The junction lattice $Z_{\tilde{X}}$ of $\tilde{X}$ is called 
\textit{junction lattice} of $X_k$ and denoted by $\tilde{L}_k$. 
Note that the parenthesis have the essential meaning
in the representation of brane configurations, namely 
$X = \left(X_{F_{u_1}}\right) \cdots \left(X_{F_{u_n}}\right)$
is a cnfluence of $\tilde{X} = X_{F_{u_1}} \cdots X_{F_{u_n}}$.

We define a mapping $v : W_1 \rightarrow \tilde{L}_k 
,\>\> J = \left[x_1,y_1;\dots;t_k;\dots;x_n,y_n\right] \mapsto v
\left(J\right)$ by 
\begin{equation*}
v\left(J\right) = \left(Q_{1,1},\dots,Q_{1,m_1};\dots;
Q_{n,1},\dots,Q_{n,m_n}\right),
\end{equation*}
where 
\begin{eqnarray}
Q_{i,j} = 
\begin{cases}
t_i &\mbox{if $u_i \in T$,} \\
\begin{pmatrix}p_{i,j} \\ q_{i,j}\end{pmatrix} \times 
K_{p_{i,j-1},q_{i,j-1}}\cdots K_{p_{i,1},q_{i,1}} 
\begin{pmatrix}x_{i} \\ y_{i}\end{pmatrix} 
\quad (1 \le j \le m_i)
&\mbox{if $u_i \notin T$.} 
\end{cases}\label{vdef}
\end{eqnarray}
for $1 \le i \le n$.
By a calculation we can show the following fact. 
\begin{equation*}
\langle J_1, J_2 \rangle = \left(v\left(J_1\right),v\left(J_2\right)
\right). 
\end{equation*}
We put $U_0 = v\left(U_1\right)$. 
Let 
\begin{equation*}
\tilde{W}_i = \left\{Q = \left(Q_{1,1},\dots,Q_{1,m_1};\dots;
Q_{n,1},\dots,Q_{n,m_n}\right) \in \bZ^{12k} \>\> \bigg | \>\> 
\sum_{j=1}^{m_i}Q_{i,j}
\begin{pmatrix}p_{i,j} \\ q_{i,j}\end{pmatrix} = 
\begin{pmatrix}0 \\ 0\end{pmatrix}\right\} 
\left(1 \le i \le n\right),
\end{equation*}
$\tilde{W}=\bigoplus_{i=1}^{n} \tilde{W}_i$, 
and $\tilde{V} = \{ J+U_0 \ \vert \ J \in \tilde{W}\}$.

\begin{Lem}\label{isomL}
{\normalfont\rmfamily
The quotient lattice $\tilde{L}_k / U_0$ is isomorphic to $L_k$.
}
\end{Lem}
\textit{Proof.} \, In Proposition $\ref{jlattice}$ we take 
$\bX_{\left[\tilde{p}_1,\tilde{q}_1\right]}\cdots
\bX_{\left[\tilde{p}_{12k},\tilde{q}_{12k}\right]} = 
X_{F_{u_1}} \cdots X_{F_{u_n}}$ instead of 
$X = \bX_{\left[p_1,q_1\right]}\cdots
\bX_{\left[p_{12k},q_{12k}\right]}$. 
It is enough to show $U = U_0$ by Proposition $\ref{jlattice}$. 
Let $\tilde{K}_i = K_{\left[\tilde{p}_i,\tilde{q}_i\right]}\>\>
\left(1 \le i \le n\right)$, and put $l = \sum_{k=1}^{i-1} m_k + j$. 
Then we obtain 
\begin{equation*}
\begin{split}
a_l &= \begin{pmatrix}p_l \\ q_l\end{pmatrix} \times 
       \tilde{K}_{l}\cdots\tilde{K}_{\left(\sum_{k=1}^{i-1} m_k\right)+1}
       \tilde{K}_{\sum_{k=1}^{i-1} m_k}\cdots\tilde{K}_{1}
       \begin{pmatrix}1\\0\end{pmatrix} \\
    &= \begin{pmatrix}p_{i,j} \\ q_{i,j}\end{pmatrix} \times 
       K_{\left[p_{i,j-1},q_{i,j-1}\right]}\cdots
       K_{\left[p_{i,1},q_{i,1}\right]}
       K_{i-1}\cdots K_1\begin{pmatrix}1\\0\end{pmatrix}
\end{split}
\end{equation*}
By Lemma $\ref{base1}$ and $\left(\ref{vdef}\right)$ we obtain $Q_a = v(J_a)$. 
Similarly we can show $Q_b = v(J_b)$. 
By $U = \langle Q_a,Q_b \rangle$ and $U_0 = \langle v\left(J_a\right), 
v\left(J_b\right)\rangle$, we have $U = U_0$. 
\textit{Q.E.D.}
\ake
\begin{Lem}\label{isomV}
{\normalfont\rmfamily
The quotient lattice $\tilde{V} / U_0$ is isomorphic to $V$.
}
\end{Lem}
\textit{Proof.} \, By the definition of $\tilde{V}$ and of $\tilde{W}$ 
we have 
\begin{equation*}
\tilde{V}/ U_0 \simeq \tilde{W} = \bigoplus_{i=1}^{n}\tilde{W}_i.
\end{equation*}
On the other side we have $\tilde{W}_i \simeq T_{u_i}\>\>
\left(1\le i \le n\right)$ by Table $\ref{table1}$ and Table $\ref{table2}$. 
Therefore we have $\tilde{V} / U_0 \simeq V$. 
\textit{Q.E.D.}
\ake
Consider the lattice embedding 
\begin{equation}
i_K : \tilde{V} / U_0 \hookrightarrow \tilde{L}_k / U_0. \label{emb2}
\end{equation}
\vskip 1mm
\begin{The}\label{ikiv}
{\normalfont\rmfamily
The lattice embedding 
$i_K : \tilde{V} / U_0 \hookrightarrow \tilde{L}_k / U_0$ is equivalent to 
the lattice embedding $i_V : V \hookrightarrow L_k$, that is, the 
following diagram commutes. 
$$\begin{CD}
i_K : \tilde{V} / U_0 @.\quad\hookrightarrow\quad @. \tilde{L}_k / U_0 \\
\quad\quad @VVV            @.                                    @VVV              \\
i_V : V               @.\quad\hookrightarrow\quad @. L_k
\end{CD}$$
}
\end{The}
\ake
\textit{Proof.} \, There exists a topological elliptic surface 
$\tilde{X_k}$ such that the brane configuration is 
$\tilde{X}=X_{F_{u_1}}\cdots X_{F_{u_n}}$. 
We denote $L_k$ of $\tilde{X_k}$ by $L'_k$. 
By definition $\tilde{L_k}$ is the junction lattice of $\tilde{X_k}$ 
and $U_0$ is the null junction $\bZ$-module of $\tilde{X_k}$. 
Let $\tilde{h} : \tilde{L_k} \rightarrow L'_k$ be the class mapping
and put $\tilde{h}(\tilde{V})=V'$. 
By Lemma \ref{Kerh}, Lemma \ref{isomL} and Lemma \ref{isomV}, 
$i_K$ is equivalent to the lattice embedding 
$i : V' \hookrightarrow L'_K$. 
On the other hand, because there exists a topological deformation 
{}from $X_k$ to $\tilde{X_k}$, $i$ is equivalent to $i_V$. 
\textit{Q.E.D.}
\ake
Let $\tilde{V}^{\perp}$ be the orthogonal complement of $\tilde{V}$ in 
$\tilde{L}_k$. 
\begin{Lem}\label{inverse}
{\normalfont\rmfamily
The mapping $v : W_1 \rightarrow \tilde{V}^{\perp}$ is isomorphic. 
}
\end{Lem}
\textit{Proof.} \, By $\left(\ref{koutensuu}\right)$ we have 
$v\left(W_1\right) \subset \tilde{V}^{\perp}$. 
Let $Q = \left(Q_{1,1},\dots,Q_{1,m_1};\dots;Q_{n,1},\dots,Q_{n,m_n}
\right)$ be an element of $\tilde{V}^{\perp}$. 
In the same way as the proof of Lemma $\ref{Kerh}$ we can show 
that $Q$ satisfies the condition $\left(\ref{vdef}\right)$. 
Hence we have $v\left(W_1\right) = \tilde{V}^{\perp}$. \\
We define the mapping $v' : \tilde{V}^{\perp} 
\rightarrow W_1 ,\>\> 
Q = \left(Q_{1,1},\dots,Q_{1,m_1};\dots;Q_{n,1},\dots,Q_{n,m_n}\right) 
\mapsto v'\left(Q\right)$ by 
\begin{equation*}
v'\left(Q\right) = \left(\sum_{j=1}^{m_1}Q_{1,j}
\left(p_{1,j},q_{1,j}\right);\dots;\sum_{j=1}^{m_n}Q_{n,j}
\left(p_{n,j},q_{n,j}\right)\right).
\end{equation*}
We have $v' \circ v = id_{W_1}$ and 
$v \circ v' = id_{\tilde{V}^{\perp}}$ by $\left(\ref{sub1}\right)$ and 
$\left(\ref{sub2}\right)$. 
We have thus proved the lemma. 
\textit{Q.E.D.}
\ake

\section{Torsion part of Mordell-Weil group, Mordell-Weil lattice, and
transcendental lattice}

\subsection{Junctions and homologies}
In this section we explain the relation between tadpole junction, 
rational tadpole junction and homology. 
\ake
First we describe the relation between a tadpole junction and homology. 

\begin{The}
{\normalfont\rmfamily
Let $W_1$ be the tadpole junction lattice of $X_k$ and $U_1$ the null 
tadpole junction \bZ-module of $X_k$. Then 
the quotient lattice ${W_1}/{U_1}$ is isomorphic to $V^{\perp}$. 
}
\end{The}
\textit{Proof.} \, By Lemma $\ref{Kerh}$, it is enough to show that the 
mapping $h$ is surjective . 
We can define the mapping $u = h \circ v^{-1} : \tilde{V}^{\perp} 
\rightarrow V^{\perp}$ by Lemma \ref{inverse}. 
It is clear that $u$ conserve the intersection number. 
Because $\mathrm{Ker} \> u = U_0$, we have $\tilde{V}^{\perp} / U_0 \simeq 
u\left(\tilde{V}^{\perp}\right)$. 
On the other side, by the meaning of $U_0$ and Theorem $\ref{ikiv}$, we have 
$\tilde{V}^{\perp} / U_0 \simeq {\left(\tilde{V} / U_0\right)}^{\perp} \simeq V^{\perp}$. 
Therefore we see that $u$ is surjective. 
It follows from this that $h$ is surjective immediately. 
\textit{Q.E.D.}
\ake
The following corollary follows from the theorem proved above and Theorem 3. 
\begin{Cor}\label{extremal}
{\normalfont\rmfamily
The transcendental lattice of extremal elliptic $K3$ surfaces is isomorphic 
to the quotient lattice ${W_1}/{U_1}$.
}
\end{Cor}
\ake
In the following we describe the relation between a rational tadpole junction 
and homology. 
\begin{Lem}
{\normalfont\rmfamily
The quotient lattice $W_2 / U_1$ is isomorphic to $L_k / V$ : 
\begin{equation}
W_2 / U_1 \simeq L_k / V . \label{gyulem}
\end{equation}
}
\end{Lem}
\textit{Proof.} \, Let the mapping $w : \tilde{L}_k 
\rightarrow W_2 ,\>\> 
Q = \left(Q_{1,1},\dots,Q_{1,m_1};\dots;Q_{n,1},\dots,Q_{n,m_n}\right) 
\mapsto w\left(Q\right)$ defined by 
\begin{equation*}
w\left(Q\right) = \left(\sum_{j=1}^{m_1}Q_{1,j}
\left(p_{1,j},q_{1,j}\right);\dots;\sum_{j=1}^{m_n}Q_{n,j}
\left(p_{n,j},q_{n,j}\right)\right)
\end{equation*}
and the mapping $\bar{w} : \tilde{L}_k / U_0 \rightarrow W_2 / U_1 $ defined 
by $J + U_0 \mapsto w\left(J\right) + U_1$. 
Let the mapping 
$\bar{v} : W_2 / U_1 \rightarrow (\tilde{L_k} \bigotimes \bQ) / 
U_0 = (\tilde{L_k} / U_0) \bigotimes \bQ$ defined by 
$J + U_1 \mapsto \tilde{v}\left(J\right) + U_1$, 
where $\tilde{v}$ is defined by the same formula $(\ref{vdef})$ for $v$. 
Let $p = \bar{v} \circ \bar{W} : \tilde{L}_k / U_0 \rightarrow 
(\tilde{L_k} / U_0) \bigotimes \bQ$. 
Then we see that $p$ is the orthogonal projection to the direction of 
$(\tilde{V} / U_0)^{\perp}$. 
By Theorem \ref{ikiv}, we have $p(\tilde{L_k} / U_0) \simeq L_k / 
V$. 
On the other hand we have 
$W_2 / U_1 \simeq p(\tilde{L_k} / U_0)$ because 
$\bar{w}(\tilde{L_k} / U_0) = \left(W_2 / U_1\right)$ and
$\bar{v} : W_2 / U_1 \rightarrow \bar{v}\left(W_2 / U_1\right)$ is 
isomorphic. 
Therefore we obtain $(\ref{gyulem})$. 
\textit{Q.E.D.}
\ake
\begin{The}
{\normalfont\rmfamily
Let $W_2$ be the rational tadpole junction lattice of $X_k$ and $U_2$ 
the rational null junction \bZ-module of $X_k$. Then 
the quotient lattice ${W_2}/{U_2}$ is isomorphic to $\left(V^{\perp}\right)^{*}$.
}
\end{The}
\textit{Proof.} \, We know the following fact \cite{OS} : 
\begin{equation*}
L_k / V \simeq \left(V^{\perp}\right)^{*} \oplus \hat{V}/V . 
\end{equation*}
The theorem follows from taking the free part of 
$\left(\ref{gyulem}\right)$. 
\textit{Q.E.D.}
\ake
It follows from the theorem proved above and Theorem 2 
\begin{Cor}\label{rational}
{\normalfont\rmfamily
The Mordell-Weil lattice of rational elliptic surfaces $X_k$ is isomorphic to 
the opposite lattice of the quotient lattice ${W_2}/{U_2}.$
}
\end{Cor}
\ake
The following theorem follows from taking the torsion part of 
$\left(\ref{gyulem}\right)$. 
\begin{The}
{\normalfont\rmfamily
Let $U_1$ be the null junction \bZ-module of $X_k$ and $U_2$ the rational null 
junction \bZ-module of $X_k$. Then 
${U_2}/{U_1}$ is isomorphic to $\hat{V}/V$.
}
\end{The}
\ake
It follows from the theorem proved above and Theorem 1 
\ake
\begin{Cor}\label{u2u1}
{\normalfont\rmfamily
The torsion part of Mordell-Weil group of elliptic surfaces is isomorphic to 
${U_2}/{U_1}$.
}
\end{Cor}

\subsection{Torsion part of Mordell-Weil group, Mordell-Weil lattice, 
transcendental lattice} \label{algo}
In this section we explain how to calculate the torsion part of 
Mordell-Weil group of elliptic surfaces, the Mordell-Weil lattice 
of rational elliptic surfaces, and the transcendental lattice and the 
Mordell-Weil group of extremal elliptic $K3$ surfaces from the given data 
of global monodromy.\\
First we describe the method to calculate the torsion part of Mordell-Weil 
group of elliptic surfaces from the given data of global monodromy. 
\begin{The}
{\normalfont\rmfamily
The torsion part of Mordell-Weil group of elliptic surfaces is 
isomorphic to $\bZ / {n_1 \bZ} \oplus \bZ / {n_2 \bZ}$. 
}
\end{The}
The theorem follows from Lemma $\ref{base1}$, Lemma $\ref{base2}$ and 
Corollary $\ref{u2u1}$, and 
the explicit algorythm to calculate the numbers $n_1$, $n_2$ is already 
given in Section \ref{null}. 
\ake
Next we describe the method to calculate the transcendental lattice of 
extremal elliptic $K3$ surfaces from the given data of global monodromy. \\
It is sufficient to get the basis of ${W_1}/{U_1}$ by Corollary 
\ref{extremal} and $\left(\ref{seki1}\right)$.
\begin{enumerate}
\item Construction of the basis $\langle J_1,J_2,\dots,J_{l-2} \rangle$ of 
$W_1$. \\
We define the $2 \times l$ matrix $R$ such that the following equation is 
equivalent to the equation $\left(\ref{pq1}\right)$.
\begin{equation*}
 R \> {}^{t}\left(x_1,y_1;\dots;t_i;\dots;x_n,y_n\right) = \begin{pmatrix} 0 \\ 0 \end{pmatrix}.
\end{equation*}
By using the elementary divisor theory we can have a 
$2 \times 2$ unimodular matrix $H$ and a $l \times l$ unimodular matrix 
$G = \left(g_1,\dots,g_{l}\right)$ such that 
\begin{equation*}
 H R G = \begin{pmatrix}u_1 & 0   & 0 & \cdots & 0 \\
                              0   & u_2 & 0 & \cdots & 0 
               \end{pmatrix}.
\end{equation*} 
Then the basis are obtained as 
\begin{equation*}
 J_1 = {}^{t}g_3, \> J_2 = {}^{t}g_4,\dots, \> J_{l-2} = {}^{t}g_{l}.
\end{equation*}
\item Construction of the basis $\langle J_a,J_b\rangle$ of $U_1$. \\
      This is given in Lemma \ref{base1}. 
\item Construction of the basis of ${W_1}/{U_1}$.\\
Let 
\begin{equation*}
 D = \begin{pmatrix} J_1 \\ \vdots \\ J_{l-2} \end{pmatrix} , \>\>
 O = \begin{pmatrix} J_a \\ J_b \end{pmatrix}.  
\end{equation*}
Then there exists a unique $2 \times (l-2)$ matrix $L$ such that 
\begin{equation*}
 L D = O.
\end{equation*}
Because the embedding $U_1 \subset W_1$ is primitive, by the elementary 
divisor theory, there exists a $2 \times 2$ unimodular matrix $P$ and a 
$(l-2) \times (l-2)$ unimodular matrix $Q$ such that 
\begin{equation*}
 P L Q = \begin{pmatrix}  1 & 0 & 0 & \cdots & 0 \\
                          0 & 1 & 0 & \cdots & 0 
         \end{pmatrix}.
\end{equation*}
Let the $(l-4) \times l$ matrix $V$ 
\begin{equation*}
 V := \begin{pmatrix} v_1 \\
                     \vdots \\
                     v_{l-4}
             \end{pmatrix} := 
 \begin{pmatrix}
 0 & 0 & 1 & 0 & \cdots & 0 \\
 0 & 0 & 0 & 1 & \cdots & 0 \\
 \vdots & \vdots & \vdots & \vdots & \ddots & \vdots \\
 0 & 0 & 0 & 0 & \cdots & 1 \\
 \end{pmatrix}
 Q^{-1}D.
\end{equation*}
Then $\langle v_1,\dots,v_{l-4} \rangle$ is the basis of ${W_1}/{U_1}$. 
\end{enumerate}
\ake
Finally we describe the method to calculate Mordell-Weil lattice of rational 
elliptic surfaces from the given data of global monodromy. \\
It is sufficient to get the basis of ${W_2}/{U_2}$ by Corollary 
\ref{rational} and $\left(\ref{seki2}\right)$.
\begin{enumerate}
\item Construction of the basis $\langle J'_1,J'_2,\dots,J'_{l-2} \rangle$ of 
$W_2$. \\
We define the $2 \times l$ matrix $R'$ such that the following equation is 
equivalent to the equation $\left(\ref{pq2}\right)$.
\begin{equation*}
 R' \> {}^{t}\left(x'_1,y'_1;\dots;t'_i;\dots;x'_n,y'_n\right) = \begin{pmatrix} 0 \\ 0 \end{pmatrix}.
\end{equation*}
By using the elementary divisor theory we can have a $2 \times 2$ 
unimodular matrix $H'$ and a $l \times l$ unimodular matrix 
$G' = \left(g'_1,\dots,g'_{l}\right)$ such that 
\begin{equation*}
 H' R' G' = \begin{pmatrix}u'_1 & 0   & 0 & \cdots & 0 \\
                                0   & u'_2 & 0 & \cdots & 0 
               \end{pmatrix}.
\end{equation*} 
Then we can take 
\begin{equation*}
 J'_1 = {}^{t}g'_3, \> J'_2 = {}^{t}g'_4,\dots, \> J'_{l-2} = {}^{t}g'_{l}.
\end{equation*}
\item Construction of the basis $\langle J'_a,J'_b\rangle$ of $U_2$ \\
This is given in Lemma \ref{base2}.
\item Construction of the basis of ${W_2}/{U_2}$.\\
Let 
\begin{equation*}
 D'= \begin{pmatrix} J'_1 \\ \vdots \\ J'_{l-2} \end{pmatrix} , \>\>
 O' = \begin{pmatrix} J'_a \\ J'_b \end{pmatrix}.  
\end{equation*}
Then there exists a unique $2 \times (l-2)$ matrix $L'$ such that 
\begin{equation*}
 L' D' = O'.
\end{equation*}
Because the embedding $U_2 \subset W_2$ is primitive, by the elementary 
divisor theory, there exists a $2 \times 2$ unimodular matrix $P'$ and 
a $(l-2) \times (l-2)$ unimodular matrix $Q'$ such that 
\begin{equation*}
 P' L' Q' = \begin{pmatrix}  1 & 0 & 0 & \cdots & 0 \\
                          0 & 1 & 0 & \cdots & 0 
         \end{pmatrix}.
\end{equation*}
Let the $(l-4) \times l$ matrix $V'$ 
\begin{equation*}
 V' := \begin{pmatrix} v'_1 \\
                     \vdots \\
                     v'_{l-4}
             \end{pmatrix} :=
 \begin{pmatrix}
 0 & 0 & 1 & 0 & \cdots & 0 \\
 0 & 0 & 0 & 1 & \cdots & 0 \\
 \vdots & \vdots & \vdots & \vdots & \ddots & \vdots \\
 0 & 0 & 0 & 0 & \cdots & 1 \\
 \end{pmatrix}
 {Q'}^{-1} D'.
\end{equation*}
Then $\langle v'_1,\dots,v'_{l-4} \rangle$ is the basis of ${W_2}/{U_2}$. 
\end{enumerate}

\begin{Exa}
{\normalfont\rmfamily
(Table \ref{config1} : No.$34$)\\
$X = \left(\bA^{4}\bB\bC\right)\left(\bA^{2}\right)\left(\bA^{2}\right)
\bB\bC$. 
\ake
$K = \left(
     \begin{pmatrix}-1 & 0 \\ 0 & -1\end{pmatrix},
     \begin{pmatrix}1 & -2 \\ 0 & 1\end{pmatrix},
     \begin{pmatrix}1 & -2 \\ 0 & 1\end{pmatrix},
     \begin{pmatrix}0 & -1 \\ 1 & 2\end{pmatrix},
     \begin{pmatrix}2 & -1 \\ 1 & 0\end{pmatrix}
     \right).$
\aki
The Mordell-Weil group $MW \cong {{A_1}^{*}}^{\oplus 2} \oplus 
\bZ / 2\bZ$ is computed as follows : 
\begin{enumerate}
\item
\begin{equation*}
\begin{split}
R' &= \begin{pmatrix} 1 & 0 & 1 & 1 & 1 & 1 \\
                        0 & 1 & 0 & 0 & -1 & 1 
     \end{pmatrix}, \\
H' &= \begin{pmatrix}1 & 0 \\ 0 & 1 \end{pmatrix}, 
G' = \begin{pmatrix} 1 & 0 & -1 & -1 & -1 & -1 \\
                     0 & 1 & 0  & 0  & 1  & -1 \\
                     0 & 0 & 1  & 0  & 0  & 0  \\
                     0 & 0 & 0  & 1  & 0  & 0  \\
                     0 & 0 & 0  & 0  & 1  & 0  \\
                     0 & 0 & 0  & 0  & 0  & 1  
     \end{pmatrix}. 
\end{split}
\end{equation*}
\item \quad $n_1 = 1, n_2 = 2$.
\begin{equation*}
N = \begin{pmatrix} 2 & 0 \\ 0 & 2 \\ 0 & 2 \\ 0 & 2 \\ 1 & 1 \\ -1 & 1 \\ 
    \end{pmatrix}, 
A = \begin{pmatrix} 0&  0& 0& 0&  1& 0\\
                    0&  1& 0& 0&  0& 0\\
                    0& -1& 1& 0&  0& 0\\
                    0& -1& 0& 1&  0& 0\\
                    1&  1& 0& 0& -2& 0\\
                    0& -1& 0& 0&  1& 1 
    \end{pmatrix}, 
B = \begin{pmatrix} 1 & -1 \\ 0 & 1\end{pmatrix}.
\end{equation*}
$J'_a = \left(2,0;0;0;-1;-1\right),J'_b = \left(-1,1;-1;-1;2;1\right)$
\item
\begin{equation*}
L' = \begin{pmatrix} 0  & 0  & -1 & -1 \\
                     -1 & -1 & 2  & 1 
     \end{pmatrix}, 
P' = \begin{pmatrix} 1 & 0 \\ 0 &1 \end{pmatrix}, 
Q' = \begin{pmatrix} -1 &  0 & -1 & -1 \\
                     0  & -1 & 2  & 1  \\
                     0  &  0 & 1  & 0  \\
                     -1 &  0 & -1 & 0 
     \end{pmatrix}.
\end{equation*}
$v'_1 = \left(-1,1,0,0,1,0\right), v'_2 = \left(0,-1,-1,0,0,1\right)$.\\
For $J' = \left(x'_1,y'_1;t'_2;t'_3;t'_4;t'_5\right)$, we have 
\begin{equation*}
\begin{split}
\langle J',J' \rangle = &-{t'}_2^2/2-{t'}_3^2/2-{t'}_4^2-{t'}_5^2
                         -x'_1{t'}_4+x'_1{t'}_5
                         -y'_1{t'}_2-y'_1{t'}_3-y'_1{t'}_4-y'_1{t'}_5 \\
                        &-{t'}_2{t'}_4+{t'}_2{t'}_5
                         -{t'}_3{t'}_4+{t'}_3{t'}_5+2{t'}_4{t'}_5. 
\end{split}
\end{equation*}
By using reduction theory, we have 
\begin{equation*}
\langle v'_1,v'_2 \rangle = \begin{pmatrix}-1 & 3/2 \\ 3/2 & -5/2\end{pmatrix} 
\cong \begin{pmatrix}-1/2 & 0 \\ 0 & -1/2\end{pmatrix} 
= \langle v'_a,v'_b \rangle, 
\end{equation*}
where $v'_a = v'_1 + v'_2 = \left(-1,0;-1;0;1;1\right), 
v'_b = -2v'_1-v'_2 = \left(2,-1;1;0;-2;-1\right)$.
\end{enumerate}
Therefore we have 
$MW \cong A_1^{\ast \oplus 2} \oplus \bZ/2\bZ$. 
}
\end{Exa}

\begin{Exa}\label{ex(a)}
{\normalfont\rmfamily
(Table \ref{config2} : No.$61\,(a)$)\\
$X = \left(\bA^{10}\right)(\bX_{[0,1]}^{6})\left(\bC^{4}\right)
(\bX_{[3,1]}^{2})\bX_{[41,11]}\bX_{[5,1]}$. 
\ake
For $J = [t_1;t_2;t_3;t_4;t_5;t_6]$, we have 
\begin{gather*}
\begin{split}
\langle J,J \rangle = &-10t_1^2-6t_2^2-4t_3^2-2t_4^2-t_5^2-t_6^2
                      +60t_1t_2+40t_1t_3+20t_1t_4+110t_1t_5+10t_1t_6\\
                      &-24t_2t_3-36t_2t_4-246t_2t_5-30t_2t_6
                      -16t_3t_4-120t_3t_5-16t_3t_6-16t_4t_5\\
                      &-4t_4t_6-14t_5t_6.
\end{split}\\
\langle v_a,v_b \rangle = \begin{pmatrix}10 & 0 \\ 0 & 12\end{pmatrix}
\quad\left(\textrm{transcendental lattice}\right),
\end{gather*}
where $v_a = [0;-7;37;-176;23;-7], v_b = [0;0;-1;4;0;-4]$. 
\ake
$U_1 = \langle J_a,J_b \rangle$, 
where $J_a = [0;-1;5;-23;3;-1], \> J_b = [1;10;-49;225;-29;5]$. 
\ake
$U_2 = \langle J'_a,J'_b \rangle$, 
where $J'_a = (0;-6;20;-46;3;-1), \> J'_b = (5;15;-48;110;-7;0)$. 
\ake
$J'_a = J_a \in U_1, \>\> J'_b \notin U_1, \>\> 2J'_b = 5J_a + J_b \in U_1,
 \>\>  MW \cong \bZ/2\bZ$. 
}
\end{Exa}

\begin{Exa}\label{ex(b)}
{\normalfont\rmfamily
(Table \ref{config2} : No.$61\,(b)$)\\
$X = \left(\bA^{10}\right)(\bB^{6})(\bX_{[4,-5]}^{4})
(\bX_{[3,-4]}^{2})\bX_{[1,-3]}\bC$. 
\ake
For $J = [t_1;t_2;t_3;t_4;t_5;t_6]$, we have 
\begin{gather*}
\begin{split}
\langle J,J \rangle = &-10t_1^2-6t_2^2-4t_3^2-2t_4^2-t_5^2-t_6^2
                      -60t_1t_2-200t_1t_3-80t_1t_4-30t_1t_5+10t_1t_6\\
                      &-24t_2t_3-12t_2t_4-12t_2t_5+12t_2t_6
                      -8t_3t_4-28t_3t_5+36t_3t_6-10t_4t_5\\
                      &+14t_4t_6+4t_5t_6.
\end{split}\\
\langle v_a,v_b \rangle = \begin{pmatrix}10 & 0 \\ 0 & 12\end{pmatrix}
\quad\left(\textrm{transcendental lattice}\right),
\end{gather*}
where $v_a = [0;3;-1;0;0;-2], v_b = [0;2;-3;6;0;0]$. 
\ake
$U_1 = \langle J_a,J_b \rangle$, 
where $J_a = [0;1;-1;2;-1;-1], \> J_b = [1;-9;8;-15;5;1]$. 
\ake
$U_2 = \langle J'_a,J'_b \rangle$, 
where $J'_a = (10;-54;32;-30;5;1), \> J'_b = (-5;24;-14;13;-2;0)$. 
\ake
$J'_a = J_b \in U_1, \>\> J'_b \notin U_1, \>\> 2J'_b = -J_a - J_b \in U_1,
 \>\>  MW \cong \bZ/2\bZ$. 
}
\end{Exa}

\begin{Rem} \label{itti}
{\normalfont\rmfamily
Let $X_1$ be the extremal semi-stable elliptic $K3$ surface of 
Example \ref{ex(a)} and $X_2$ be the one of Example \ref{ex(b)}. 
Then the triplet $\left(V, V^{\perp}, \hat{V}/V\right)$ of $X_1$ coincides 
with the one of $X_2$, however the embeddings $i_V$ of these 
surface are not equivalent. 
The difference of these embeddings can be seen by looking at the 
structure of sections. \\
Let $X$ be an extremal semi-stable elliptic $K3$ surface which have $I_{n_1}
,\dots,I_{n_6}$ type fibres, and $\psi$ be the homomorphism from $MW(X)$ to 
$\bZ/n_1\bZ \times \cdots \times \bZ/n_6\bZ$ given in \cite{MP}, i.e., 
$\psi\left(s\right) = \left(a_1,\dots,a_6\right)$, where $a_i$ specify 
the irreducible component that the section $s$ intersects the 
corresponding singular fibre. 
Let $s_1$ be the generator of $MW(X_1)$, and then we have 
$\psi\left(s_1\right) = \left(5,3,0,0,0,0\right)$ \cite{MP}\cite{ABTZ}. 
We can consider that $J'_b = \left(5;15;-48;110;-7;0\right)$ corresponds 
to $s_1$ at this case. In fact we have 
$J'_b \equiv \psi\left(s_1\right)\>\> \textrm{mod} \> \left(10,6,4,2,1,1\right)$. 
Similarly for the surface of $X_2$, the generator of $MW(X_2)$ $s_2$ 
satisfy $\psi\left(s_2\right) = \left(5,0,2,1,0,0\right)$ and 
$J'_b = \left(-5;24;-14;13;-2;0\right)$. 
Then we also have 
$J'_b \equiv \psi\left(s_1\right)\>\> \textrm{mod} \> \left(10,6,4,2,1,1\right)$. 
}
\end{Rem}

\section{The global monodromies of rational elliptic surfaces and of 
extremal elliptic $K3$ surfaces}

\subsection{Some working hypothesises for classifications of global monodromies}
In this section we discuss the global monodromies 
$K = \left(K_1,\dots,K_n\right)$
of elliptic surfaces. 
\aki
Let 
\begin{equation*}
M_1 = \left\{K \mid K \>\> \textrm{satifies the condition (a)\,} 
      \right\}, \>\>\>
M_2 = \left\{K \mid K \>\> \textrm{satifies the condition (b)\,} 
      \right\}, 
\end{equation*}
where the conditions (a), (b) are as follows : 
\begin{list}{}{}
\item[(a)]There exists a \textit{complex} elliptic surface such that the 
global monodromy is $K$ and the Euler number is $12k$. 
\item[(b)]There exists a \textit{topological} elliptic surface such that 
the global monodromy is $K$ and the Euler number is $12k$. 
\end{list}
Let 
\begin{gather*}
E_1 = \left\{i_V : V \hookrightarrow L_k \mid i_V \>\> \textrm{satifies 
      the condition (\setcounter{enumiii}{1}\theenumiii)} 
      \right\}, \\
E_2 = \left\{i_V : V \hookrightarrow L_k \mid i_V \>\> \textrm{satifies 
      the condition (\setcounter{enumiii}{2}\theenumiii)} 
      \right\}, 
\end{gather*}
where the conditions (\setcounter{enumiii}{1}\theenumiii), 
(\setcounter{enumiii}{2}\theenumiii) is as follows : 
\begin{list}{}{}
\item[(\setcounter{enumiii}{1}\theenumiii)]
There exists a \textit{complex} elliptic surface which realize the 
embedding $i_V$. 
\item[(\setcounter{enumiii}{2}\theenumiii)]
There exists a \textit{topological} elliptic surface which realize the 
embedding $i_V$. 
\end{list}
Let the mapping $F : M_2 \rightarrow E_2 $ defined by 
$K \mapsto i_K : V \hookrightarrow L_k$. 
\ake
We define the following three kinds of operations \setcounter{part}{1}\thepart, 
\setcounter{part}{2}\thepart, \setcounter{part}{3}\thepart\, acting on 
monodromy data $K = \left(K_1,\dots,K_n\right) \in M$. 
\begin{list}{}{}
\item[\setcounter{part}{1} \thepart.]
The elementary transformation : \\
$R_i : K \longmapsto R_i\left(K\right) \quad (1 \le i \le n-1)$. 
\item[\setcounter{part}{2} \thepart.]
The action of $P \in SL\left(2,\bZ\right)$ : \\
$K \longmapsto K'=\left(K'_1,\dots,K'_n\right)$, where 
$K'_i = P K_i P^{-1}\left(1 \le i \le n\right)$. 
\item[\setcounter{part}{3} \thepart.]
$K \longmapsto K'=\left(K'_1,\dots,K'_{n+1}\right)$. \\
$\left(K_1,\ldots,K_{i-1},K_C{K_A}^{h},K_{i+1},\ldots,K_n\right) 
\longmapsto \left(K_1,\ldots,K_{i-1},{K_A}^{h},K_C,K_{i+1},\ldots,K_n\right)$, 
where $h = 2,3,4$. 
\end{list}
If $K$ transfer to $K'$ by a finite sequence of the 
operations \setcounter{part}{1}\thepart, 
\setcounter{part}{2}\thepart, \setcounter{part}{3}\thepart\ and of the 
inverse operations of them, then we denote $K \approx K'$. 
It is easy to see that $F\left(K\right) = F\left(K'\right)$ if 
$K \approx K'$. We anticipate that the converse is true. 
\ake
\textbf{Hypothesis 1.} \\
$(\mathrm{H_1})$
Let $K$, $K'$ be the elements of $M_2$. 
Then we have $K \approx K'$ if 
$F\left(K\right) = F\left(K'\right)$. 
\ake
Let 
\begin{equation*}
\tilde{M_2} = \left\{K \mid K \>\> \textrm{satifies the condition (d)\,} 
      \right\} \subseteq M_2, 
\end{equation*}
where the condition (c) is as follows : 
\begin{list}{}{}
\item[(c)]there exists a topological elliptic surface such that the global 
           monodromy is $K$ and the Euler number is $12k$, 
           and which has no singular fibre of type 
           \setcounter{part}{2}$\thepart$, 
           \setcounter{part}{3}$\thepart$, 
           \setcounter{part}{4}$\thepart$. 
\end{list}
It follows from Hypothesis 1. immediately.
\ake
\textbf{Hypothesis 2.}\,
$(\mathrm{H_2})$
The mapping 
$F|_{\tilde{M_2}} : \tilde{M_2} \rightarrow E_2$ is bijective. 
\ake
\textbf{Hypothesis 3.}\,
$(\mathrm{H_3})$
The mapping 
$F|_{M_1 \cap \tilde{M_2}} : M_1 \cap \tilde{M_2} \rightarrow E_1$ is surjective. 
\ake
Let 
\begin{equation*}
M = M_1 \cap \tilde{M_2}, \quad M' = F^{-1}\left(E_1\right) \cap \tilde{M_2}. 
\end{equation*}
We assume Hypothesis $1.$ and Hypothesis $2.$.
\ake
\indent We discuss the classification of $K \in M_1$. 
It follows from $(\mathrm{H_2})$ and $(\mathrm{H_3})$ 
that the mapping 
$F|_M : M \rightarrow E_1$ is a bijection. 
On the other hand we see that the mapping 
$F|_{M'} : M' \rightarrow E_1$ is a bijection by $(\mathrm{H_2})$. 
Therefore, under $(\mathrm{H_2})$ and $(\mathrm{H_3})$, we conclude 
\begin{equation}
M = M'. \label{conclu1}
\end{equation}
Let $V$ be the direct sum of some root lattices of type $A$,$D$,$E$. \\
Let 
\begin{equation*}
D = \left\{V \mid V \>\>\textrm{satisfies the 
    condition} \> \textrm{(\setcounter{enumiii}{3}\theenumiii)} 
    \right\}, 
\end{equation*}
where 
\begin{list}{}{}
\item[(\setcounter{enumiii}{3}\theenumiii)]
there exists a complex elliptic surface such that the trivial lattice is 
$V \oplus U'$. 
\end{list}
For $V \in D$, we define the equation $e_{V}$ as follows : 
\begin{equation*}
e_{V} : K_n \cdots K_1 = I, 
\end{equation*}
where each $K_i$ is conjugate to the local monodromy matrix of the 
singular fibre of type $I_i$, $I_i^{*}$,
 \setcounter{part}{2}$\thepart^{*}$, \setcounter{part}{3}$\thepart^{*}$, 
\setcounter{part}{4}$\thepart^{*}$ corresponding to the lattice $V$. 
By the definition of $M'$, we have 
\begin{equation*}
M' = \left\{K = \left(K_1,\dots,K_n\right) \mid K \>\> \textrm{is a 
     solution of} \>\> e_V, V \in D
    \right\}. \label{conclu3}
\end{equation*}
\ake
Let 
\begin{equation*}
R = \left\{(V, V^{\perp}, \hat{V}/V) \mid (V, V^{\perp}, \hat{V}/V) 
       \>\> \textrm{satifies the condition} \>\> 
      (\mathrm{\setcounter{enumiii}{4}\theenumiii})
      \right\}, 
\end{equation*}
where 
\begin{list}{}{}
\item[(\setcounter{enumiii}{4}\theenumiii)]
there exists a complex elliptic surface which realize 
$(V, V^{\perp}, \hat{V}/V)$. 
\end{list}
In the case of rational elliptic surafaces ($k = 1$) and of extremal 
elliptic $K3$ surfaces (special case of $k = 2$)\, $R$ \, is 
Oguiso-Shioda's table and Shimada-Zhang's table respectively. 
Our problem is to obtain the list of global monodromies which 
reproduce Oguiso-Shioda's table and Shimada-Zhang's table. 
\subsection{A classification table}
For the problem raised at the end of previous subsction, we obtained 
a (tentative) result, which is summerized in Table \ref{config1} and 
Table \ref{config2}. \\
The problem is to find the configuration which reproduce the 
Oguiso-Shioda's table and Shimada-Zhang's table. 
Let us explain how we solved this problem. 
\ake
Some of the configuratons are obatined simply 
by the confluence process of known solutions 
$\bA^{8}\bB\bC\bB\bC$ and $\bA^{8}\bB\bC\bB\bC\bA^{8}\bB\bC\bB\bC$ 
(Table  \ref{config1}, No.$1$ - $8$ etc). 
To obtain the other configurations we use the elementary transformations. 
For example, the configuration of No.34 in Table 3, the transformation 
is given as follows : 
\begin{equation*}
\begin{split}
 & \bA^8\bB\bC\bB\bC \sim \bA^4\bA^3\bA\bB\bC\bB\bC 
 \sim \bA^4\bA^3\bB{\bX_{[0,1]}}\bC\bB\bC \sim \bA^4\bA^3\bB\bC\bA\bB\bC 
 \sim \bA^4\bA^2\bB\bC\bA^2\bB\bC \\
 & \sim \bA^4\bB\bC\bA^2\bA^2\bB\bC
\end{split}
\end{equation*}
We obtained almost all the configurations in the tables by this method. 
For each configuration, we have calculated the pair 
$(V^{\perp}, \hat{V}/V)$  by the method of section 4.2. 
For the cases such as the case No.18 in Table \ref{config2}, 
we need to find the configuration corresponding for each data of 
$(V^{\perp}, \hat{V}/V)$. 
In these cases, one of the configurations is sometimes hard to obtain and 
we need a lot of try and error. 
The systematic method of the caluclation of the data 
$(V^{\perp}, \hat{V}/V)$ was useful to complete this procedure. 
\ake
In the following we consider the meanings of Table \ref{config1}. and of Table 
\ref{config2}. under the above hypothesises. \\
Let 
\begin{equation*}
S_1 = \left\{K \mid K \>\> \textrm{appears in Table} \>\> 
\ref{config1}. 
      \right\}, \quad
S_2 = \left\{K \mid K \>\> \textrm{appears in Table} \>\> 
\ref{config2}.
      \right\}.
\end{equation*}
\indent First we consider the case of rational elliptic surfaces. 
We assume $\mathrm{(H_2)}$. 
One can show that $\mathrm{(H_3)}$ is true in this case \cite{UP}. 
On the other hand it is known that the classification of 
$i : V \hookrightarrow L_1 \in E_1$ coincide with the classification of 
$(V, V^{\perp}, \hat{V}/V) \in R_1$\, (= Oguiso-shioda' table) \cite{OS}. 
Hence we have $M' = S_1$. 
Therefore, by $\left(\ref{conclu1}\right)$, under $\mathrm{(H_2)}$ we conclude 
\begin{equation*}
M = S_1.
\end{equation*}
\indent In the following we consider the case of extremal elliptic $K3$ 
surfaces. 
We assume $(\mathrm{H_2})$. 
In this case there are no possibility of the difference between $M$ and 
$M_1$ except for No.$297$ \cite{Y}. 
Hence we see that $\mathrm{(H_3)}$ is true except for No.$297$. 
Therefore we have $M_1 = M'$ except for No.$297$. 
On the other hand it turns out that the classification of 
$(V, V^{\perp}, \hat{V}/V) \in R_1$ is slightly diffrent from 
the classification of $i : V \hookrightarrow L_2 \in E_1$. 
See Remark \ref{itti}. 
\ake
\textbf{Assumption.}\,
$(\mathrm{A})$
There are no difference between the classification of $
(V, V^{\perp}, \hat{V}/V) \in R_1$ and the one of 
$i : V \hookrightarrow L_2 \in E_1$ except for No.$61$. 
\ake
Under $\mathrm{(H_2)}$ and $\mathrm{(A)}$ we conclude 
\begin{equation*}
M_1 = S_2 \quad \textrm{except for\,No.}297 \>\> \textrm{and \,No.}61. 
\end{equation*}

In the Table 4 the data $(a, b, c)$ in the last columun repesent the
transcidental lattice, namely
$$
matrix...
$$
 
The trivial lattice $V$ is recovered from the brane configurations
by the following rule:
$$
(X^n) -> A_{n-1}
$$
The numbering is the same as that of Simada-Zhang. 

\vskip6mm\noindent
{\bf Acknowledgements}

\vskip2mm
We would like to thank Y. Yamada for many valuable discussions and 
useful conservations. 

\renewcommand{\arraystretch}{1.3}
\begin{table}
\begin{center}
\caption{Global monodromies of rational elliptic surfaces}\label{config1}
\aki

\end{center}
\end{table}
\clearpage
\begin{table}
\begin{center}
\addtocounter{table}{-1}
\caption{Global monodromies of rational elliptic surfaces}
\aki
%
\end{center}
\end{table}
\clearpage
\begin{table}[t]
\begin{center}
\addtocounter{table}{-1}
\caption{Global monodromies of rational elliptic surfaces}
\aki
%
\end{equation*}

\clearpage

\begin{table}[htbp]
 \begin{center}
 \caption{Global monodromies of extremal elliptic $K3$ surfaces}\label{config2}
 \aki
  \begin{tabular}{|l|c|c|c|}
    \hline
      No & global monodromy & $MW$ & $a\quad\quad b\quad\quad c$ \\
    \hline\hline
     $1$ & $(\bA^{4})(\bX_{[3,-2]}^{4})(\bX_{[10,-7]}^{4})
     (\bX_{[7,-5]}^{4})(\bX_{[4,-3]}^{4})(\bB^{4})$ & $\bZ/4\bZ \times 
     \bZ/4\bZ$  & $4\quad\quad 0\quad\quad 4$  \\  
    \hline\hline
     $2$ & $(\bA^{5})(\bB^{5})(\bX_{[2,-3]}^{5})(\bX_{[3,-5]}^{5})
     (\bX_{[0,1]}^{2})(\bX_{[10,3]}^{2})$ & $\bZ/5\bZ$ & $10\quad\quad 
     0 \quad\quad10$ \\  
    \hline\hline
     $3$ & $(\bA^{5})(\bX_{[3,-2]}^{5})(\bX_{[13,-9]}^{4})
     (\bX_{[10,-7]}^{4})(\bX_{[7,-5]}^{3})(\bB^{3})$ & $\left(0\right)$ 
     & $60\quad\quad 0\quad\quad60$ \\  
    \hline\hline
     $4$ & $(\bA^{6})(\bB^{6})(\bX_{[4,-5]}^{6})(\bX_{[3,-4]}^{2})
     (\bX_{[3,-5]}^{2})(\bX_{[0,1]}^{2})$ & $\bZ/2\bZ \times \bZ/6\bZ$ & 
     $2\quad\quad 0\quad\quad6$ \\  
    \hline\hline
     $5$ & $(\bA^{6})(\bX_{[1,2]}^{6})(\bX_{[6,11]}^{3})
     (\bX_{[4,7]}^{3})(\bX_{[3,5]}^{3})(\bC^{3})$ & $\bZ/3\bZ \times 
     \bZ/3\bZ$ &  $6\quad\quad 0\quad\quad6$  \\  
    \hline\hline
     $6$ & $(\bA^{6})(\bB^{6})(\bX_{[82,-99]}^{6})(\bX_{[4,-5]}^{4})
     \bX_{[29,-37]}\bX_{[5,-7]}$ & $\bZ/6\bZ$ & $4\quad\quad 0\quad\quad6$ \\  
    \hline\hline
     $7$ & $(\bA^{6})(\bB^{6})(\bX_{[4,-5]}^{4})(\bX_{[0,1]}^{4})
     (\bX_{[3,8]}^{2})(\bC^{2})$ & $\bZ/2\bZ \times \bZ/2\bZ$ & 
     $12\quad\quad  0\quad\quad12$ \\
\hline
  \end{tabular}
 \end{center}
\end{table}
\clearpage
\begin{table}[htbp]
 \begin{center}
 \addtocounter{table}{-1}
 \caption{Global monodromies of extremal elliptic $K3$ surfaces}
 \aki
  \begin{tabular}{|l|c|c|c|}
    \hline
      No & global monodromy & $MW$ & $a\quad\quad b\quad\quad c$ \\
    \hline\hline
     $24$ & $(\bA^{7})(\bX_{[0,1]}^{6})(\bX_{[5,7]}^{5})(\bC^{4})
     \bX_{[7,2]}\bX_{[20,3]}$ & $\left(0\right)$ & $12\quad\quad 0 
     \quad\quad70$  \\  
    \hline\hline
     $25$ & $(\bA^{8})(\bX_{[0,1]}^{8})(\bX_{[1,7]}^{2})
     (\bX_{[1,5]}^{2})(\bX_{[1,3]}^{2})(\bC^{2})$ & $\bZ/2\bZ \times 
     \bZ/4\bZ$ & $4\quad\quad 0 \quad\quad4$ \\  
    \hline\hline
     $26.\ (1)$ & $(\bA^{8})(\bB^{8})(\bX_{[0,1]}^{3})(\bC^{3})
     \bX_{[7,3]}\bX_{[4,1]}$ & $\left(0\right)$ & $24\quad\quad 0 
     \quad\quad24$  \\  
    \hline
     $\phantom{26.}\ (2)$ & $(\bA^{8})(\bX_{[0,1]}^{8})(\bC^{3})
     (\bX_{[5,3]}^{3})\bX_{[19,11]}\bX_{[3,1]}$ & $\bZ/2\bZ$ & 
     $12\quad\quad 0 \quad\quad12$ \\  
    \hline\hline
     $27$ & $(\bA^{8})(\bX_{[2,-1]}^{8})(\bB^{4})(\bX_{[5,1]}^{2})
     \bX_{[39,7]}\bX_{[29,5]}$ & $\bZ/8\bZ$ & $2\quad\quad 0 \quad\quad4$ \\  
    \hline\hline
     $28$ & $(\bA^{8})(\bX_{[2,-1]}^{4})(\bX_{[5,-3]}^{4})
     (\bX_{[8,-5]}^{4})(\bB^{2})(\bX_{[3,1]}^{2})$ & $\bZ/2\bZ \times 
     \bZ/4\bZ$ & $4\quad\quad 0 \quad\quad8$ \\  
    \hline\hline
     $29$ & $(\bA^{8})(\bB^{4})(\bX_{[1,-2]}^{4})(\bX_{[2,-5]}^{4})
     (\bX_{[1,-3]}^{3})\bC$ & $\bZ/4\bZ$ & $4\quad\quad 0 \quad\quad24$ \\
\hline
  \end{tabular}
 \end{center}
\end{table}
\clearpage
\begin{table}[htbp]
 \begin{center}
 \addtocounter{table}{-1}
 \caption{Global monodromies of extremal elliptic $K3$ surfaces}
 \aki
  \begin{tabular}{|l|c|c|c|}
    \hline
      No & global monodromy & $MW$ & $a\quad\quad b\quad\quad c$ \\
    \hline\hline
     $45$ & $(\bA^{9})(\bX_{[2,-1]}^{5})(\bX_{[7,-4]}^{3})
     (\bX_{[5,-3]}^{3})(\bB^{3})\bC$ & $\bZ/3\bZ$ & $12\quad\quad 3 
     \quad\quad12$ \\  
    \hline\hline
     $46$ & $(\bA^{9})(\bX_{[0,1]}^{5})(\bC^{4})(\bX_{[3,2]}^{3})
     (\bX_{[7,4]}^{2})\bX_{[3,1]}$ & $\left(0\right)$ & $6\quad\quad 0 
     \quad\quad180$ \\  
    \hline\hline
     $47$ & $(\bA^{9})(\bB^{6})(\bX_{[0,1]}^{3})(\bC^{3})
     (\bX_{[5,2]}^{2})\bX_{[4,1]}$ & $\bZ/3\bZ$ & $6\quad\quad 0 
     \quad\quad18$  \\  
    \hline\hline
     $48$ & $(\bA^{9})(\bB^{6})(\bX_{[0,1]}^{4})(\bC^{3})\bX_{[12,7]}
     \bX_{[3,1]}$ & $\bZ/3\bZ$ & $4\quad\quad 0 \quad\quad18$ \\  
    \hline\hline
     $49$ & $(\bA^{9})(\bB^{6})(\bX_{[0,1]}^{5})(\bX_{[1,2]}^{2})
     \bX_{[2,1]}\bX_{[15,2]}$ & $\left(0\right)$ & $18\quad\quad 0 
     \quad\quad30$  \\  
    \hline\hline
     $50$ & $(\bA^{9})(\bB^{7})(\bX_{[3,-5]}^{3})(\bX_{[0,1]}^{2})
     (\bX_{[4,1]}^{2})\bX_{[15,2]}$ & $\left(0\right)$ & $18\quad\quad 0 
     \quad\quad42$ \\  
    \hline\hline
     $51$ & $(\bA^{9})(\bB^{7})(\bX_{[0,1]}^{4})(\bC^{2})
     \bX_{[7,2]}\bX_{[6,1]}$ & $\left(0\right)$ & $10\quad\quad 4 
     \quad\quad52$  \\
\hline
  \end{tabular}
 \end{center}
\end{table}
\clearpage
\begin{table}[htbp]
 \begin{center}
 \addtocounter{table}{-1}
 \caption{Global monodromies of extremal elliptic $K3$ surfaces}
 \aki
  \begin{tabular}{|l|c|c|c|}
    \hline
      No & global monodromy & $MW$ & $a\quad\quad b\quad\quad c$ \\
    \hline\hline
     $67$ & $(\bA^{10})(\bB^{9})(\bC^{2})\bX_{[32,19]}\bX_{[12,7]}
       \bX_{[3,1]}$ & $\left(0\right)$ & $10\quad\quad 0 \quad\quad18$ \\
    \hline\hline
     $68$ & $(\bA^{11})(\bB^{4})(\bX_{[1,-3]}^{4})(\bX_{[3,-10]}^{3})
     \bX_{[1,-5]}\bC$ & $\left(0\right)$ & $24\quad\quad 12 \quad\quad28$ \\
    \hline\hline
     $69$ & $(\bA^{11})(\bB^{4})(\bX_{[5,-7]}^{3})(\bX_{[2,-3]}^{3})
     (\bC^{2})\bX_{[3,1]}$ & $\left(0\right)$ & $12\quad\quad 0 
     \quad\quad66$ \\
    \hline\hline
     $70$ & $(\bA^{11})(\bB^{5})(\bX_{[15,-19]}^{5})\bX_{[2,-3]}
     \bX_{[2,1]}\bX_{[11,2]}$ & $\left(0\right)$ & $10\quad\quad 5 
     \quad\quad30$  \\
    \hline\hline
     $71.\ (1)$ & $(\bA^{11})(\bB^{5})(\bX_{[31,-39]}^{3})(\bC^{3})
     \bX_{[17,12]}\bX_{[8,5]}$ & $\left(0\right)$ & $6\quad\quad 3 
     \quad\quad84$  \\  
    \hline
     $\phantom{71.}\ (2)$ &
     $(\bA^{11})(\bB^{5})(\bX_{[2,-3]}^{3})(\bX_{[3,-5]}^{3})
     \bX_{[1,-3]}\bC$ & $\left(0\right)$ & $24\quad\quad 9 \quad\quad24$ \\
    \hline\hline
     $72$ & $(\bA^{11})(\bB^{5})(\bX_{[1,3]}^{3})(\bX_{[2,5]}^{2})
     (\bC^{2})\bX_{[19,2]}$ & $\left(0\right)$ & $2\quad\quad 0 
     \quad\quad330$  \\
\hline
  \end{tabular}
 \end{center}
\end{table}
\clearpage
\begin{table}[htbp]
 \begin{center}
 \addtocounter{table}{-1}
 \caption{Global monodromies of extremal elliptic $K3$ surfaces}
 \aki
  \begin{tabular}{|l|c|c|c|}
    \hline
      No & global monodromy & $MW$ & $a\quad\quad b\quad\quad c$ \\
    \hline\hline
     $85$ & $(\bA^{12})(\bB^{5})(\bX_{[3,-4]}^{2})(\bX_{[3,-5]}^{2})
     (\bX_{[0,1]}^{2})\bX_{[5,1]}$ & $\bZ/2\bZ$ & $6\quad\quad 0 
     \quad\quad20$  \\  
    \hline\hline
     $86$ & $(\bA^{12})(\bC^{5})(\bX_{[5,4]}^{3})(\bX_{[3,2]}^{2})
     \bX_{[3,1]}\bX_{[10,1]}$ & $\left(0\right)$ & $12\quad\quad 0 
     \quad\quad30$  \\  
    \hline\hline
     $87.\ (1)$ & $(\bA^{12})(\bB^{6})(\bX_{[4,-5]}^{2})(\bX_{[2,-3]}^{2})
     \bX_{[1,-3]}\bC$ & $\bZ/2\bZ$ & $6\quad\quad 0 \quad\quad12$ \\  
    \hline
     $\phantom{87.}\ (2)$ & $(\bA^{12})(\bX_{[2,1]}^{6})(\bX_{[9,4]}^{2})
     (\bX_{[3,1]}^{2})\bX_{[24,5]}\bX_{[6,1]}$ & $\bZ/6\bZ$ & 
     $2\quad\quad 0 \quad\quad4$ \\  
    \hline\hline
     $88$ & $(\bA^{12})(\bX_{[3,1]}^{6})(\bX_{[4,1]}^{3})\bX_{[23,5]}
     \bX_{[26,5]}\bX_{[11,2]}$ & $\bZ/3\bZ$ & $4\quad\quad 0 \quad\quad6$ \\  
    \hline\hline
     $89$ & $(\bA^{12})(\bX_{[3,-1]}^{7})(\bX_{[7,1]}^{2})\bX_{[53,7]}
     \bX_{[61,8]}\bX_{[17,2]}$ & $\left(0\right)$ & $4\quad\quad 0 
     \quad\quad42$ \\  
    \hline\hline
     $90$ & $(\bA^{13})(\bX_{[5,-1]}^{3})(\bX_{[3,1]}^{3})
     (\bX_{[4,1]}^{2})(\bX_{[37,6]}^{2})\bX_{[13,2]}$ & $\left(0\right)$ 
     & $12\quad\quad 6 \quad\quad42$ \\
\hline
  \end{tabular}
 \end{center}
\end{table}
\clearpage
\begin{table}[htbp]
 \begin{center}
 \addtocounter{table}{-1}
 \caption{Global monodromies of extremal elliptic $K3$ surfaces}
 \aki
  \begin{tabular}{|l|c|c|c|}
    \hline
      No & global monodromy & $MW$ & $a\quad\quad b\quad\quad c$ \\
    \hline\hline
     $105$ & $(\bA^{15})(\bC^{5})\bX_{[23,19]}\bX_{[28,23]}\bX_{[16,13]}
       \bX_{[3,2]}$ & $\left(0\right)$ & $10\quad\quad 5 \quad\quad10$ \\  
    \hline\hline
     $106$ & $(\bA^{16})(\bB^{2})(\bC^{2})(\bX_{[3,1]}^{2})
     \bX_{[5,1]}\bX_{[13,1]}$ & $\bZ/4\bZ$ & $2\quad\quad 0 \quad\quad4$ \\  
    \hline\hline
     $107.\ (1)$ & $(\bA^{16})(\bX_{[0,1]}^{3})
       (\bC^{2})\bX_{[7,2]}\bX_{[6,1]}\bX_{[9,1]}$ & $\left(0\right)$ & 
       $10\quad\quad 2 \quad\quad10$ \\  
    \hline
     $\phantom{107.}\ (2)$ & $(\bA^{16})(\bB^{3})(\bC^{2})
     \bX_{[19,11]}\bX_{[9,5]}\bX_{[3,1]}$ & $\bZ/2\bZ$ & $4\quad\quad 0 
     \quad\quad6$ \\  
    \hline\hline
     $108$ & $(\bA^{16})(\bB^{4})\bX_{[5,-7]}\bX_{[3,-5]}\bX_{[1,-3]}
       \bC$ & $\bZ/4\bZ$ & $2\quad\quad 0 \quad\quad2$ \\  
    \hline\hline
     $109.\ (1)$ & $(\bA^{17})(\bC^{2})(\bX_{[3,1]}^{2})\bX_{[19,5]}
     \bX_{[21,5]}\bX_{[9,2]}$ & $\left(0\right)$ & $2\quad\quad 0 
     \quad\quad34$ \\  
    \hline
     $\phantom{109.}\ (2)$ & $(\bA^{17})(\bX_{[3,1]}^{2})(\bX_{[21,5]}^{2})
     \bX_{[9,2]}\bX_{[8,1]}\bX_{[18,1]}$ & $\left(0\right)$ & $4\quad\quad 
     2 \quad\quad18$ \\
\hline
  \end{tabular}
 \end{center}
\end{table}
\clearpage
\begin{table}[htbp]
 \begin{center}
 \addtocounter{table}{-1}
 \caption{Global monodromies of extremal elliptic $K3$ surfaces}
 \aki
  \begin{tabular}{|l|c|c|c|}
    \hline
      No & global monodromy & $MW$ & $a\quad\quad b\quad\quad c$ \\
    \hline\hline
     $126$ & $(\bA^{5}\bB\bC)(\bX_{[3,-2]}^{9})(\bB^{6})
     \bX_{[5,-7]}\bX_{[3,-5]}$ & $\left(0\right)$ & $12\quad\quad 0 
     \quad\quad18$ \\  
    \hline\hline
     $127$ & $(\bA^{5}\bB\bC)(\bB^{10})(\bX_{[7,-8]}^{3})
     (\bX_{[6,-7]}^{3})\bX_{[4,-5]}$ & $\left(0\right)$ & $6\quad\quad 0 
     \quad\quad60$ \\  
    \hline\hline
     $128$ & $(\bA^{5}\bB\bC)(\bB^{10})(\bX_{[3,-4]}^{3})
     (\bX_{[5,-7]}^{2})(\bX_{[1,-2]}^{2})$ & $\bZ/2\bZ$ & $2\quad\quad 0 
     \quad\quad60$ \\  
    \hline\hline
     $129$ & $(\bA^{5}\bB\bC)(\bB^{10})(\bX_{[7,-8]}^{4})
     (\bX_{[13,-15]}^{2})\bX_{[5,-6]}$ & $\bZ/2\bZ$ & $8\quad\quad 4 
     \quad\quad12$ \\  
    \hline\hline
     $130$ & $(\bA^{5}\bB\bC)(\bX_{[3,-2]}^{10})(\bB^{5})
     \bX_{[5,-8]}\bX_{[4,-7]}$ & $\left(0\right)$ & $10\quad\quad 0 
     \quad\quad20$ \\  
    \hline\hline
     $131$ & $(\bA^{5}\bB\bC)(\bX_{[0,1]}^{11})(\bX_{[1,8]}^{3})
     (\bX_{[1,7]}^{2})\bX_{[1,3]}$ & $\left(0\right)$ & $14\quad\quad 4 
     \quad\quad20$ \\  
    \hline\hline
     $132$ & $(\bA^{5}\bB\bC)(\bB^{12})(\bX_{[9,-10]}^{2})
     (\bX_{[7,-8]}^{2})\bX_{[5,-6]}$ & $\bZ/4\bZ$ & $2\quad\quad 0 
     \quad\quad6$  \\
\hline
  \end{tabular}
 \end{center}
\end{table}
\clearpage
\begin{table}[htbp]
 \begin{center}
 \addtocounter{table}{-1}
 \caption{Global monodromies of extremal elliptic $K3$ surfaces}
 \aki
  \begin{tabular}{|l|c|c|c|}
    \hline
      No & global monodromy & $MW$ & $a\quad\quad b\quad\quad c$ \\
    \hline\hline
     $149$ & $(\bA^{6}\bB\bC)(\bX_{[2,-1]}^{10})(\bB^{4})
     \bX_{[3,-5]}\bX_{[1,-3]}$ & $\bZ/2\bZ$ & $4\quad\quad 0 \quad\quad10$ \\  
    \hline\hline
     $150$ & $(\bA^{6}\bB\bC)(\bX_{[2,-1]}^{11})(\bB^{3})
     \bX_{[1,-4]}\bX_{[1,-7]}$ & $\left(0\right)$ & $6\quad\quad 0 
     \quad\quad22$  \\  
    \hline\hline
     $151$ & $(\bA^{6}\bB\bC)(\bB^{12})(\bX_{[1,-2]}^{2})
     \bX_{[1,-4]}\bX_{[1,2]}$ & $\bZ/2\bZ$ & $4\quad\quad 0 \quad\quad6$ \\  
    \hline\hline
     $152$ & $(\bA^{6}\bB\bC)(\bB^{13})\bX_{[10,-11]}
     \bX_{[5,-6]}\bX_{[2,-3]}$ & $\left(0\right)$ & $4\quad\quad 2 
     \quad\quad14$  \\  
    \hline\hline
     $153$ & $(\bA^{6}\bB\bC)(\bX_{[1,-3]}^{5}\bX_{[0,1]}\bX_{[2,-5]})
     (\bX_{[0,1]}^{6})(\bC^{3})$ & $\bZ/2\bZ$ & $6\quad\quad 0 \quad\quad12$ \\  
    \hline\hline
     $154$ & $(\bA^{6}\bB\bC)(\bX_{[1,-2]}^{5}\bX_{[0,1]}\bX_{[2,-3]})
     (\bX_{[0,1]}^{8})\bX_{[1,2]}$ & $\bZ/2\bZ$ & $4\quad\quad 0 
     \quad\quad8$ \\   
    \hline\hline
     $155$ & $(\bA^{7}\bB\bC)(\bB^{7}\bX_{[1,-2]}\bA)(\bX_{[1,-2]}^{3})
     (\bX_{[0,1]}^{3})$ & $\left(0\right)$ & $12\quad\quad 0 \quad\quad12$ \\
\hline  
  \end{tabular}
 \end{center}
\end{table}
\clearpage
\begin{table}[htbp]
 \begin{center}
 \addtocounter{table}{-1}
 \caption{Global monodromies of extremal elliptic $K3$ surfaces}
 \aki
  \begin{tabular}{|l|c|c|c|}
    \hline
      No & global monodromy & $MW$ & $a\quad\quad b\quad\quad c$ \\
    \hline\hline
     $174$ & $(\bA^{8}\bB\bC)(\bB^{7})(\bX_{[4,-5]}^{3})(\bX_{[7,-9]}^{3})
     \bX_{[1,-2]}$ & $\left(0\right)$ & $12\quad\quad 6 \quad\quad24$ \\  
    \hline\hline
     $175$ & $(\bA^{8}\bB\bC)(\bB^{8})(\bX_{[2,-3]}^{3})(\bX_{[0,1]}^{2})
     \bX_{[8,3]}$ & $\bZ/2\bZ$ & $2\quad\quad 0 \quad\quad24$ \\  
    \hline\hline
     $176$ & $(\bA^{8}\bB\bC)(\bB^{10})(\bX_{[0,1]}^{2})\bX_{[7,6]}
     \bX_{[3,2]}$ & $\bZ/2\bZ$ & $2\quad\quad 0 \quad\quad10$ \\  
    \hline\hline
     $177$ & $(\bA^{8}\bB\bC)(\bX_{[3,1]}^{5}\bA\bX_{[5,2]})
     (\bX_{[5,1]}^{5}\bA\bX_{[9,2]})$ & $\bZ/2\bZ$ & $4\quad\quad 0 
     \quad\quad4$  \\  
    \hline\hline
     $178$ & $(\bA^{8}\bB\bC)(\bB^{6}\bX_{[0,1]}\bX_{[2,-1})(\bX_{[0,1]}^{4})
     (\bC^{2})$ & $\bZ/2\bZ \times \bZ/2\bZ$ & $2\quad\quad 0 \quad\quad4$ \\  
    \hline\hline
     $179$ & $(\bA^{9}\bB\bC)(\bX_{[4,1]}^{9}\bX_{[5,1]}\bX_{[3,1]})
     \bX_{[9,2]}\bX_{[7,1]}$ & $\left(0\right)$ & $4\quad\quad 0 
     \quad\quad4$ \\   
    \hline\hline
     $180$ & $(\bA^{9}\bB\bC)(\bB^{5})(\bX_{[7,-9]}^{3})(\bX_{[3,-4]}^{3})
     (\bX_{[0,1]}^{2})$ & $\left(0\right)$ & $12\quad\quad 0 \quad\quad30$ \\
\hline  
  \end{tabular}
 \end{center}
\end{table}
\clearpage
\begin{table}[htbp]
 \begin{center}
 \addtocounter{table}{-1}
 \caption{Global monodromies of extremal elliptic $K3$ surfaces}
 \aki
  \begin{tabular}{|l|c|c|c|}
    \hline
      No & global monodromy & $MW$ & $a\quad\quad b\quad\quad c$ \\
    \hline\hline
     $199$ & 
     $(\bA^{11}\bB\bC)(\bX_{[3,-1]}^{5})(\bB^{3})(\bX_{[0,1]}^{2})
     \bX_{[7,2]}$ & $\left(0\right)$ & $6 \quad\quad 0 \quad\quad20$ \\  
    \hline\hline
     $200$ & $(\bA^{11}\bB\bC)(\bX_{[2,-1]}^{6})(\bB^{3})\bX_{[2,-5]}\bC$ & 
     $\left(0\right)$ & $6\quad\quad 0 \quad\quad12$ \\  
    \hline\hline
     $201$ & 
     $(\bA^{11}\bB\bC)(\bB^{7})(\bX_{[1,-2]}^{2})\bX_{[2,-7]}\bC$ & 
     $\left(0\right)$ & $6\quad\quad 2 \quad\quad10$ \\  
    \hline\hline
     $202$ & 
     $(\bA^{12}\bB\bC)(\bX_{[3,-1]}^{3})(\bC^{3})(\bX_{[3,2]}^{2})
     (\bX_{[4,1]}^{2})$ & $\bZ/2\bZ$ & $6\quad\quad 0 \quad\quad6$ \\  
    \hline\hline
     $203$ & 
     $(\bA^{12}\bB\bC)(\bX_{[2,-1]}^{4})(\bB^{3})(\bX_{[1,-3]}^{2})\bC$ 
     & $\bZ/2\bZ$ & $4\quad\quad 0 \quad\quad6$ \\  
    \hline\hline
     $204$ & 
     $(\bA^{12}\bB\bC)(\bB^{5})(\bX_{[3,-4]}^{2})(\bX_{[1,-2]}^{2})\bC$ 
     & $\bZ/2\bZ$ & $4\quad\quad 2 \quad\quad6$ \\  
    \hline\hline
     $205$ & 
     $(\bA^{12}\bB\bC)(\bB^{5}\bX_{[2,-3]}\bX_{[0,1]})(\bC^{2})
     \bX_{[3,1]}$ & $\bZ/2\bZ$ & $2\quad\quad 0 \quad\quad4$ \\
\hline  
  \end{tabular}
 \end{center}
\end{table}
\clearpage
\begin{table}[htbp]
 \begin{center}
 \addtocounter{table}{-1}
 \caption{Global monodromies of extremal elliptic $K3$ surfaces}
 \aki
  \begin{tabular}{|l|c|c|c|}
    \hline
      No & global monodromy & $MW$ & $a\quad\quad b\quad\quad c$ \\
    \hline\hline
     $224$ & 
     $(\bA^{5}\bB\bC^{2})(\bA^{6})(\bX_{[3,-1]}^{4})(\bB^{3})
     (\bX_{[1,-2]}^{3})$ & $\bZ/3\bZ$ & $6\quad\quad 0 \quad\quad12$ \\  
    \hline\hline
     $225$ & 
     $(\bA^{5}\bB\bC^{2})(\bX_{[8,3]}^{6})(\bA^{5})(\bX_{[2,-1]}^{4})
     \bX_{[2,-3]}$ & $\left(0\right)$ & $12\quad\quad 0 \quad\quad30$ \\  
    \hline\hline
     $226$ & 
     $(\bA^{5}\bB\bC^{2})(\bA^{5}\bX_{[3,-1]}\bB^{2})(\bX_{[0,1]}^{7})
     \bX_{[1,2]}$ & $\left(0\right)$ & $6\quad\quad 3 \quad\quad12$ \\  
    \hline\hline
     $227.\ (1)$ & 
     $(\bA^{5}\bB\bC^{2})(\bA^{7})(\bX_{[9,-2]}^{4})(\bX_{[4,-1]}^{3})
     (\bX_{[2,-1]}^{2})$ & $\left(0\right)$ & $6\quad\quad 0 \quad\quad84$ \\  
    \hline\hline
     $\phantom{227.}\ (2)$ & 
     $(\bA^{5}\bB\bC^{2})(\bX_{[3,1]}^{7})(\bA^{4})(\bX_{[1,-2]}^{3})
     (\bX_{[1,-3]}^{2})$ & $\left(0\right)$ & $12\quad\quad 0 \quad\quad42$ \\  
    \hline\hline
     $228$ & 
     $(\bA^{5}\bB\bC^{2})(\bA^{7})(\bX_{[4,-1]}^{5})(\bB^{2})
     (\bX_{[1,-3]}^{2})$ & $\left(0\right)$ & $20\quad\quad 10 \quad\quad26$ \\  
    \hline\hline
     $229$ & 
     $(\bA^{5}\bB\bC^{2})(\bX_{[3,1]}^{7})(\bA^{5})(\bX_{[9,-5]}^{3})
     \bX_{[5,-3]}$ & $\left(0\right)$ & $18 \quad\quad 3 \quad\quad18$ \\
\hline  
  \end{tabular}
 \end{center}
\end{table}
\clearpage
\begin{table}[htbp]
 \begin{center}
 \addtocounter{table}{-1}
 \caption{Global monodromies of extremal elliptic $K3$ surfaces}
 \aki
  \begin{tabular}{|l|c|c|c|}
    \hline
      No & global monodromy & $MW$ & $a\quad\quad b\quad\quad c$ \\
    \hline\hline
     $247$ & 
     $(\bA^{5}\bB\bC^{2})(\bX_{[4,1]}^{6}\bA\bX_{[7,2]})(\bA^{5})
     (\bX_{[0,1]}^{3})$ & $\left(0\right)$ & $6\quad\quad 0 \quad\quad30$ \\  
    \hline\hline
     $248$ & 
     $(\bA^{5}\bB\bC^{2})(\bX_{[3,1]}^{6}\bA\bX_{[5,2]})(\bA^{7})
     \bX_{[2,-1]}$ & $\left(0\right)$ & $4\quad\quad 2\quad\quad22$ \\  
    \hline\hline
     $249$ & 
     $(\bA^{5}\bB\bC^{2})(\bX_{[3,1]}^{7}\bX_{[4,1]}\bX_{[2,1]})(\bA^{5})
     (\bB^{2})$ & $\left(0\right)$ & $4\quad\quad 0 \quad\quad30$ \\  
    \hline\hline
     $250$ & 
     $(\bA^{5}\bB\bC^{2})(\bA^{7}\bB\bC)(\bX_{[0,1]}^{5}\bA\bX_{[1,-2]})$ & 
     $\left(0\right)$ & $4\quad\quad 0 \quad\quad12$ \\  
    \hline\hline
     $251$ & 
     $(\bA^{5}\bB\bC^{2})(\bA^{8}\bX_{[3,1]}\bX_{[5,1]})(\bX_{[0,1]}^{5})
     \bX_{[1,3]}$ & $\left(0\right)$ & $8\quad\quad 2 \quad\quad8$ \\  
    \hline\hline
     $252$ & 
     $(\bA^{5}\bB\bC^{2})(\bA^{9}\bX_{[2,1]}\bX_{[4,1]})(\bX_{[2,-1]}^{3})
     (\bB^{2})$ & $\left(0\right)$ & $6\quad\quad 0 \quad\quad12$ \\  
    \hline\hline
     $253$ & 
     $(\bA^{5}\bB\bC^{2})(\bX_{[3,1]}^{9}\bX_{[10,3]}\bX_{[4,1]})
     (\bA^{4})\bX_{[1,-2]}$ & $\left(0\right)$ & $4\quad\quad 0 \quad\quad12$ 
     \\
\hline  
  \end{tabular}
 \end{center}
\end{table}
\clearpage
\begin{table}[htbp]
 \begin{center}
 \addtocounter{table}{-1}
 \caption{Global monodromies of extremal elliptic $K3$ surfaces}
 \aki
  \begin{tabular}{|l|c|c|c|}
    \hline
      No & global monodromy & $MW$ & $a\quad\quad b\quad\quad c$ \\
    \hline\hline
     $271$ & 
     $(\bA^{6}\bB\bC^{2})(\bX_{[2,1]}^{8})(\bX_{[3,1]}^{4})(\bX_{[7,2]}^{2})
     \bX_{[5,1]}$ & $\bZ/2\bZ$ & $4\quad\quad 0 \quad\quad8$ \\  
    \hline\hline
     $272$ & 
     $(\bA^{6}\bB\bC^{2})(\bX_{[7,4]}^{8})(\bX_{[9,5]}^{5})\bX_{[20,11]}
     \bX_{[4,1]}$ & $\left(0\right)$ & $6\quad\quad 2 \quad\quad14$ \\  
    \hline\hline
     $273$ & 
     $(\bA^{6}\bB\bC^{2})(\bX_{[3,1]}^{9})(\bA^{3})(\bX_{[1,2]}^{2})
     \bX_{[3,4]}$ & $\left(0\right)$ & $6\quad\quad 0 \quad\quad18$ \\  
    \hline\hline
     $274$ & 
     $(\bA^{6}\bB\bC^{2})(\bX_{[3,1]}^{9})(\bA^{4})\bX_{[1,-2]}\bX_{[1,2]}$ & 
     $\left(0\right)$ & $4\quad\quad 0 \quad\quad18$ \\  
    \hline\hline
     $275$ & 
     $(\bA^{6}\bB\bC^{2})(\bA^{10})(\bX_{[7,-1]}^{2})(\bB^{2})
     \bX_{[1,-3]}$ & $\bZ/2\bZ$ & $2\quad\quad 0 \quad\quad10$ \\  
    \hline\hline
     $276.\ (1)$ & 
     $(\bA^{6}\bB\bC^{2})(\bX_{[3,1]}^{10})(\bA^{3})\bX_{[1,3]}\bX_{[2,3]}$ & 
     $\left(0\right)$ & $6\quad\quad 0 \quad\quad10$ \\  
    \hline\hline
     $\phantom{276.}\ (2)$ & 
     $(\bA^{6}\bB\bC^{2})(\bX_{[2,1]}^{10})(\bX_{[3,1]}^{3})\bX_{[11,3]}
     \bX_{[5,1]}$ & $\bZ/2\bZ$ & $4\quad\quad 1 \quad\quad4$ \\
\hline  
  \end{tabular}
 \end{center}
\end{table}
\clearpage
\begin{table}[htbp]
 \begin{center}
 \addtocounter{table}{-1}
 \caption{Global monodromies of extremal elliptic $K3$ surfaces}
 \aki
  \begin{tabular}{|l|c|c|c|}
    \hline
      No & global monodromy & $MW$ & $a\quad\quad b\quad\quad c$ \\
    \hline\hline
     $294$ & 
     $(\bA^{6}\bB\bC^{2})(\bA^{5}\bX_{[2,-1]}\bX_{[0,1]}^{2})(\bC^{6})
     \bX_{[3,2]}$ & $\left(0\right)$ & $6\quad\quad 0 \quad\quad6$ \\  
    \hline\hline
     $295$ & 
     $(\bA^{6}\bB\bC^{2})(\bA^{5}\bX_{[0,1]}\bX_{[2,1]}^{2})
     (\bA^{5}\bB\bC)$ & $\left(0\right)$ & $2\quad\quad 0 \quad\quad12$ \\  
    \hline\hline
     $296$ & 
     $(\bA^{7}\bB\bC^{2})(\bC^{7}\bX_{[5,4]}\bX_{[3,2]}^{2})(\bX_{[3,1]}^{2})
     (\bX_{[4,1]}^{2})$ & $\left(0\right)$ & $2\quad\quad 0 \quad\quad2$ \\  
    \hline\hline
     $297.\ (\alpha)$ & 
     $(\bA^{7}\bB\bC^{2})(\bA^{7}\bB\bC^{2})(\bA^{3})\bC$ & 
     $\left(0\right)$ & $2\quad\quad 1 \quad\quad2$ \\  
    \hline\hline
     $\phantom{297.}\ (\beta)$ & 
     $(\bA^{7}\bB\bC^{2})(\bA^{7}\bB\bC^{2})(\bA^{3}\bC)$ & 
     $\left(0\right)$ & $2\quad\quad 1 \quad\quad2$ \\  
    \hline\hline
     $298$ & 
     $(\bA^{7}\bB\bC^{2})(\bA^{4})(\bX_{[2,1]}^{4})(\bX_{[5,2]}^{3})
     (\bX_{[8,3]}^{3})$ & $\left(0\right)$ & $12\quad\quad 0 \quad\quad12$ \\  
    \hline\hline
     $299$ & 
     $(\bA^{7}\bB\bC^{2})(\bA^{5})(\bX_{[2,1]}^{5})(\bX_{[5,2]}^{2})
     (\bX_{[11,4]}^{2})$ & $\left(0\right)$ & 10$\quad\quad 0 \quad\quad10$ \\
\hline  
  \end{tabular}
 \end{center}
\end{table}
\clearpage
\begin{table}[htbp]
 \begin{center}
 \addtocounter{table}{-1}
 \caption{Global monodromies of extremal elliptic $K3$ surfaces}
 \aki
  \begin{tabular}{|l|c|c|c|}
    \hline
      No & global monodromy & $MW$ & $a\quad\quad b\quad\quad c$ \\
    \hline\hline
     $319$ & 
     $(\bA^{7}\bB\bC^{2})(\bA^{9}\bB\bC)(\bX_{[2,-1]}^{2})\bX_{[0,1]}$ & 
     $\left(0\right)$ & $2\quad\quad 0 \quad\quad4$ \\  
    \hline\hline
     $320$ & 
     $(\bA^{7}\bB\bC^{2})(\bA^{10}\bB\bC)\bX_{[7,-2]}\bX_{[2,-1]}$ & 
     $\left(0\right)$ & $2\quad\quad 0 \quad\quad2$ \\  
    \hline\hline
     $321$ & 
     $(\bA^{7}\bB\bC^{2})(\bA^{5}\bX_{[2,1]}\bX_{[4,1]}^{2})(\bA^{4})
     (\bX_{[2,1]}^{2})$ & $\left(0\right)$ & $2\quad\quad 0 \quad\quad12$ \\  
    \hline\hline
     $322$ & 
     $(\bA^{7}\bB\bC^{2})(\bA^{5}\bC\bX_{[3,1]}^{2})(\bA^{5})\bC$ & 
     $\left(0\right)$ & $2\quad\quad 1 \quad\quad8$ \\  
    \hline\hline
     $323$ & 
     $(\bA^{7}\bB\bC^{2})(\bA^{5}\bX_{[0,1]}\bX_{[2,1]}^{2})
     (\bA^{4}\bX_{[0,1]}\bX_{[2,1]})$ & $\left(0\right)$ & $4\quad\quad 
     2 \quad\quad4$ \\  
    \hline\hline
     $324$ & 
     $(\bA^{7}\bB\bC^{2})(\bA^{6}\bC\bX_{[3,1]}^{2})(\bA^{3})
     (\bX_{[2,1]}^{2})$ & $\left(0\right)$ & $2\quad\quad 0 \quad\quad6$ \\  
    \hline\hline
     $325$ & 
     $(\bA^{7}\bB\bC^{2})(\bA^{6}\bX_{[0,1]}\bX_{[2,1]}^{2})
     (\bA^{4})\bC$ & $\left(0\right)$ & $2\quad\quad 0 \quad\quad4$ \\  
    \hline  
  \end{tabular}
 \end{center}
\end{table}


\end{document}